\documentclass[11pt]{article}
\usepackage{amsmath}
\usepackage{latexsym,amssymb,amsmath,amsfonts,amsthm}
\usepackage{graphicx}
\usepackage{float}

\makeatletter
\newif\ifmsbmloaded@
\def\loadmsbm{\msbmloaded@true
  \font\tenmsb=msbm10 scaled 1\@ptsize00
  \font\sevenmsb=msbm7 scaled 1\@ptsize00
  \font\fivemsb=msbm5 scaled 1\@ptsize00
  \alloc@8\fam\chardef\sixt@@n\msbfam
  \textfont\msbfam=\tenmsb
  \scriptfont\msbfam=\sevenmsb
  \scriptscriptfont\msbfam=\fivemsb
  }
\def\nonmatherr@#1{\errmessage%
{LateX error: \string#1\space allowed only in math mode}}
\def\Bbb{\relax\ifmmode\expandafter\Bbb@\else
  \expandafter\nonmatherr@\expandafter\Bbb\fi}
\def\Bbb@#1{{\Bbb@@{#1}}}
\def\Bbb@@#1{\fam\msbfam\relax#1}
\@addtoreset{equation}{section}

\makeatother \loadmsbm

\def\R{\Bbb R}

\def\N{\Bbb N}

\def\C{\Bbb C}

\def\no{\noindent}
\def\f{\frac}

\def\vep{\varepsilon}

\def\p{\partial}

\def\endproof{\hphantom{MM}\hfill\llap{$\square$}\goodbreak}

\def\s{\vskip0.25cm \noindent}
\newcommand{\beq}{\begin{equation}}
\newcommand{\eeq}{\end{equation}}
\newcommand{\ben}{\begin{eqnarray}}
\newcommand{\een}{\end{eqnarray}}
\newcommand{\beno}{\begin{eqnarray*}}
\newcommand{\eeno}{\end{eqnarray*}}

\newtheorem{Theorem}{Theorem}[section]
\newtheorem{Lemma}[Theorem]{Lemma}

\newtheorem{Proposition}{Proposition}[section]
\newtheorem{Corollary}{Corollary}[section]
\begin{document}
\title{Orbital stability of traveling waves for the one-dimensional Gross-Pitaevskii equation}

\author{Patrick G\'{e}rard$^{\dag}$ and Zhifei Zhang $^{\dag\ddag}$\\
{\small $ ^\dag$  Universit\'{e} Paris-Sud, Math\'{e}matiques, B\^{a}t. 425, 91405 Orsay Cedex,
%}\\{\small
France.} \\
{\small E-mail: patrick.gerard@math.u-psud.fr}\\[2mm]
{\small $ ^\ddag$ School of  Mathematical Science, Peking University, Beijing 100871, China}\\
{\small E-mail: zfzhang@math.pku.edu.cn}\\[2mm]\\
}

\maketitle

\begin{abstract}
In this paper, we prove the nonlinear orbital stability of the stationary traveling wave of  the one-dimensional Gross-Pitaevskii equation
by using Zakharov-Shabat's inverse scattering method.
\end{abstract}

\section{Introduction}

The Gross-Pitaevskii equation

\begin{equation} \label{1.1}
\left\{ \begin{aligned}
         &i u_t+\Delta u=(|u|^2-1)u, \quad x\in \R^d \\
         &u(0,x)=u_0(x),
                          \end{aligned} \right.
\end{equation}
which models the dynamics of Bose-Einstein condensates, superfluids,
has received a lot of interest during the recent years. For a recent
state of the art, we refer to the Proceedings  \cite{FS} and to
references therein. At least formally, (\ref{1.1}) can be seen as a
Hamiltonian evolution equation associated to the Ginzburg-Landau
energy \beq H(u)=\int_{\R^d}\f12 |\nabla _xu|^2+\f14(|u|^2-1)^2dx,
\label{1.2} \eeq defined on the energy space \beno E=\{u\in
H^1_{\textrm{loc}}(\R^d): \nabla _x u\in L^2(\R^d), |u|^2-1\in
L^2(\R^d)\}. \eeno

The unusual conditions at infinity imposed by the finiteness of the
Ginzburg-Landau energy give rise to the existence of many traveling
waves solutions to (\ref{1.1}). In what follows, we shall restrict
our purpose to space dimension $d=1$, refering to \cite{Bet1},
\cite{Bet2}, \cite{Bet3} and references therein for the case $d=2$
and $d=3$. In the one-dimensional case, solutions of the form
$$
u(t,x)=U(x-ct)
$$
are completely characterized by solving an ordinary differential
equation (see e.g. \cite{Bet3}). Besides constant solutions of
modulus 1, which correspond to the solutions with $H(u)=0$, these
solutions are given by \beq U_c(x)=\sqrt{1-\f {c^2}
2}\tanh\Bigl(\sqrt{1-\f {c^2} 2}\f x{\sqrt{2}}\Bigr)+i\f c
{\sqrt{2}} \ ,\ \label{1.3} \eeq up to a multiplicative constant of
modulus $1$. An important question is then the stability of such
objects for a natural distance on the energy space $E$, say
\begin{equation}\label{distance}
d_E(u,v)=|u(0)-v(0)|+\| u'-v'\| _{L^2(\R )}+\| |u|^2-|v|^2 \| _{L^2(\R )}\ .
\end{equation}
First of all, let us mention that the notion of stability in this case has to be properly defined,
taking into account the existence of a continuum of traveling waves corresponding to a continuum of velocities.
For example, using the formula (\ref{1.3}), it is clear that $d_E(U_c,U_{c_0})$ tends to 0 as $c$ tends to $c_0$,
but if $c\neq c_0$, we have
$$
\lim_{t\rightarrow +\infty}\int_{\R}|U'_c(x-ct)-U_{c_0}(x-c_0t)|^2dx=\|U_c'\|_{L^2}^2+\|U_{c_0}'\|_{L^2}^2.
$$
For this reason, we shall say that $U_c$ is orbitally stable for the distance $d_E$ on $E$ if, denoting by $\tau_y$ the translation operator
$$
\tau_yf(x)=f(x-y),
$$
we have, for every solution $u$ of (\ref{1.1})
$$
\sup _{t\in \R}\, \inf_{y\in \R}\, d_E(\tau_{y}u(t),U_c)\rightarrow
0, \quad \textrm{as}\quad d_E(u(0),U_c)\rightarrow 0.
$$
In the case of the velocity $c\neq 0$, Lin \cite{Lin} proved the
orbital stability of $U_c$ for (\ref{1.1}) by using the
Grillakis-Shatah-Strauss theory. His proof is based on the
hydrodynamical form of (\ref{1.1}): the solution is written as
$u=(1-r)^\f12e^{i\theta}$ so that the equations expressed in terms
of new variables $(r,\theta_x)$ turn into a Hamiltonian system which
fits in the framework of \cite{GSS}. Then Lin reduces the orbital
stability of $U_c$ to the condition
$$
\f d {dc}P(U_c)<0,
$$
where the renormalized momentum  $P(u)$ is defined by
$$
P(u)=\int_{\R}\textrm{Im}(\overline{u}u')\Bigl(1-\f1 {|u|^2}\Bigr)dx.
$$
This approach is valid for non-zero velocities $c$, because in this
case $U_c$ does not vanish on $\R$. A major difficulty in extending
the above approach to the case of zero velocity is that $U_0(x)$
vanishes at some point so that $P(U_0)$ is not defined. In
\cite{DG}, Di Menza and Gallo proved the linear stability of $U_0$
under $H^1$ perturbations, and performed numerical computations
which suggest the nonlinear orbital stability. A very recent result
by B\'ethuel, Gravejat, Saut and Smets \cite{Bet4} proved a weak
form of orbital stability of $U_0$  for $d_E$, allowing to
renormalize the solution by factors of modulus $1$. In the present
paper, we prove that this renormalization by factors of modulus $1$
is useless, at least for sufficiently  smooth and decaying
perturbations. Our main result is the following  nonlinear orbital
stability of $U_0$ for (\ref{1.1}).

\begin{Theorem}\label{main}
Assume that the initial datum of (\ref{1.1}) has the form
$$
u_0(x)=U_0(x)+\varepsilon u_1(x),\qquad U_0(x)=\tanh(\f x {\sqrt{2}}),
$$
where $u_1(x)$ satisfies the following condition \beq \sup _{x\in
\R}|<x>^4\p^k u_1(x)|\le 1,\qquad \textrm{for} \quad k\le 3
.\label{1.4} \eeq Then if $\varepsilon>0$ is small enough, there
exists a unique solution $u(t,x)$ of (\ref{1.1}) such that \ben
\forall t\in \R\ ,\ \exists y(t)\in \R\ ,\
\|\tau_{y(t)}u(t,.)-U_0\|_{L^\infty}\le C\vep,\qquad
\textrm{for}\quad 0\le t<+\infty. \een
\end{Theorem}
Using a functional analytic argument from \cite{Bet4}, Theorem
\ref{main} easily yields the orbital stability for $d_E$, at least
for sufficiently smooth and decaying perturbations.
\begin{Corollary}\label{cor}
For every $\delta >0$, there exists $\vep >0$ such that, if
$$ \sup _{x\in
\R}|<x>^4\p^k u_1(x)|\le 1,\qquad \textrm{for} \quad k\le 3 ,$$ then
the solution $u$ of (\ref{1.1}) satisfies
$$\forall t\in \R\ ,\ \exists y(t)\in \R\ ,\ d_E(\tau _{y(t)}u(t),U_0)\, \leq \delta .$$
\end{Corollary}
More precisely, we shall prove that, in Theorem \ref {main} and
Corollary \ref{cor}, one can choose  $y(t)=\beta t$ and we shall
give an interpretation of the parameter $\beta $. Our strategy for
proving Theorem \ref{main} follows the inverse scattering method as
developed by  Zakharov and Shabat in \cite{Zak}. Recall that this
method  is based on the following observation. Let $u$ be a function
of $(t,x)$.  Denote by $u^*$ the complex conjugate of $u$. Set
 \ben \label{2.2} L_u=
i\left(\begin{array}{cc}
1+\sqrt{3} & 0\\
0&  1-\sqrt{3}
\end{array}\right)\f \p {\p x}+\left(\begin{array}{cc}
0 \quad& u^*\\
u \quad& 0
\end{array}\right)\ ,
\een and
\ben \label{4.1} B_u= -\sqrt{3}\left(\begin{array}{cc}
1\quad & 0\\
0\quad & 1
\end{array}\right)\f {\p^2} {\p x^2}+\left(\begin{array}{cc}
\f{|u|^2-1}{\sqrt{3}+1} & iu^*_x\\
-i u_x & \f{|u|^2-1}{\sqrt{3}-1}
\end{array}\right)\ .
\een It is easy to verify that \beno [L_u,B_u]=
\left(\begin{array}{cc}
0 & -u_{xx}^*+(|u|^2-1)u^*\\
u_{xx}-(|u|^2-1)u &  0
\end{array}\right).
\eeno Thus $u=u(t)$ satisfies the Gross-Pitaevskii equation
(\ref{1.1}) if and only if it satisfies the operator evolution
equation \ben\label{4.2} \f d {dt}L_u=i[L_u, B_u]. \een The above
evolution equation implies that spectral properties of $L_u$ are
easily handled as $t$ varies. Then the solution $u$ at time $t$ is
recovered from spectral data of $L_{u(t)}$ from a system of integral
equations. We shall follow this procedure step by step in the
perturbation context of Theorem \ref {main} and deduce the parameter
$\beta $ from the spectral data of $L_{u_0}$. \vskip0.25cm \noindent
This paper is organized as follows. In section 2, we classically
discuss the properties of generalized eigenfunctions  for $L_u$, if
$u-U_0$ is sufficiently smooth and decaying at infinity, and we
introduce a representation formula for these Jost solutions
involving a kernel $\Psi $ which will be the center of our analysis.
In section 3 we discuss the properties of transition (or scattering)
coefficients for $L_u$, in particular in the perturbation context of
Theorem \ref {main}, while section 4 is devoted to the evolution of
these coefficients deduced from equation (\ref{4.2}). Section 5 is
devoted to establishing the fundamental system of Marchenko
equations which gives the kernel $\Psi $ from the transition
coefficients. In section 6, we use the system of Marchenko equations
to prove further information about transition coefficients, in
particular the fact that the transmission coefficient admits exactly
one zero $\lambda _0$ if $u$ is a smooth and decaying perturbation
of $U_0$. Finally, Theorem \ref{main} is proved in section 7, where
it is shown that the translation vector $y(t)$ at time $t$ in can be
taken as $y(t)=-2\lambda _0t$. Appendix A is devoted to the proof of
two technical lemmas, while appendix B explains how to derive
Corollary \ref{cor} from Theorem \ref{main} and a compactness
argument in \cite{Bet4}. \s Throughout this paper, $z^*$ denotes the
complex conjugate of the complex number $z$.

\s {\bf Acknowledgements.} The first author would like to thank
J.-C. Saut for suggesting this problem. This paper was written while
the second author was visiting the Mathematics Department at Orsay
as a postdoctoral fellow. He would like to thank the hospitality and
support of the Department. The second author is partially supported
by NSF of China under Grant 10601002.

\section{Jost solutions and their properties}

In this section, we assume that $u$ is a $C^4$ function on $\R $
satisfying \beq\label{estu} \sup _{x\geq 0}\left
|\frac{d^k}{dx^k}\left (u(x)-1\right )\right |\langle x\rangle^4
+\sup _{x\leq 0}\left |\frac{d^k}{dx^k}\left (u(x)+1\right )\right
|\langle x\rangle ^4<+\infty\ ,\ 0\leq k\leq 3. \eeq In view of the
Lax pair framework recalled at the end of the introduction, the
scattering problem associated with the Gross-Pitaevskii equation is
\ben\label{2.1} L_u\chi=E\chi,\quad \chi=\left(\begin{array}{c}
\chi_1\\
\chi_2
\end{array}\right),\quad E\in \R \ ,
\een where $L_u$ is defined by (\ref{2.2}). We make the change of
variables \beno \chi_1=(\sqrt{3}-1)^\f12e^{i\f {Ex} 2}v_1,\quad
\chi_2=(\sqrt{3}+1)^\f12e^{i\f {Ex} 2}v_2. \eeno Then (\ref{2.1}) is
reduced to
\begin{align} \label{2.3} \left\{
\begin{aligned}
i\f {\p v_1} {\p x}+q^*v_2=\lambda v_1,\\
-i\f {\p v_2} {\p x}+qv_1=\lambda v_2,
\end{aligned}
\right. \end{align}
where $\lambda=\f {\sqrt{3}} 2E,\,\, q=\f {\sqrt{2}} 2u.$
Introducing the matrices
\begin{eqnarray*}
M=\left (\begin{array}{cc}1&0\\0&-1\end{array}\right )\  ,\  Q(x)=\left (\begin{array}{cc}0&q^*(x)\\ q(x)&0\end{array}\right )\ ,
\end{eqnarray*}
system (\ref{2.3}) reads
\begin{equation}\label{matrix2.3}
iM\partial _xv+Qv-\lambda v=0\ .
\end{equation}
Notice that system (\ref{2.3}) is invariant with  respect to
the involution
$$
v=\left(\begin{array}{c}
v_1\\
v_2
\end{array}\right)\rightarrow \widetilde{v}=\left(\begin{array}{c}
                                              v_2^*\\
                                               v_1^*
                                              \end{array}\right).
$$
In other words, $\widetilde{v}$ is also a solution of (\ref{2.3}) if $v$ is a solution of (\ref{2.3}).
Moreover, if both $v$ and $w$ are solutions of (\ref{2.3}), then their Wronskian defined by
$$
\{v,w\}=v_1w_2-v_2w_1
$$
does not depend on $x$.

By (\ref{estu}), we find that
$q\rightarrow \f {\sqrt{2}} 2$ as $x\rightarrow +\infty$,
and $q\rightarrow -\f {\sqrt{2}} 2$ as $x\rightarrow -\infty$.
Set
\ben\label{2.4}
X_1^+=e^{-i\zeta x}\left(\begin{array}{c}
1\\
\sqrt{2}(\lambda-\zeta)
\end{array}\right), \quad
X_2^+=e^{i\zeta x}\left(\begin{array}{c}
\sqrt{2}(\lambda-\zeta)\\
1
\end{array}\right).
\een which are solutions of (\ref{2.3}) with $q=\f {\sqrt{2}} 2$.
Here $\zeta \in \R $ and satisfies $\lambda ^2-\zeta ^2=\f 12$.
Hence $X_1^+, X_2^+$ and all the functions we are going to define
are strictly speaking functions on the hyperbola \ben \label{H}
H=\{
(\lambda ,\zeta )\in \R ^2\ |\ \lambda ^2-\zeta ^2=\f 12\} . \een
Notice that $\zeta $ is a coordinate on each of the branches $H_{\pm
}=H\cap \{ \pm \lambda >0\} $ of $H$.  In what follows, we will
often make the traditional abuse of notation which consists in
suppressing  the dependence on $\zeta$, so that these functions
appear as double-valued functions of $\lambda $. Similarly, we set
\ben\label{2.5} X_1^-=e^{-i\zeta x}\left(\begin{array}{c}
1\\
-\sqrt{2}(\lambda-\zeta)
\end{array}\right), \quad
X_2^-=e^{i\zeta x}\left(\begin{array}{c}
-\sqrt{2}(\lambda-\zeta)\\
1
\end{array}\right),
\een which are solutions of (\ref{2.3}) with $q=-\f {\sqrt{2}} 2$.
\vskip 0.25cm \no
 The Jost solutions $\psi_1$ and $\psi_2$ are the
solutions of the system (\ref{2.3}) with the following asymptotic
forms at infinity
$$
\psi_1\sim X_1^+,\qquad \psi_2\sim X_2^+, \quad \textrm{as} \quad x\rightarrow+\infty.
$$
Similarly,
$$
\varphi_1\sim  X_1^-,\qquad \varphi_2\sim  X_2^-, \quad \textrm{as} \quad x\rightarrow-\infty.
$$
Let us recall why these functions are well defined. We decompose $Q(x)$ as
$$Q(x)=Q^++R^+(x)\ ,\ Q^+=\left (\begin{array}{cc}0&\frac{\sqrt 2}{2}\\\frac{\sqrt 2}{2}&0\end{array}\right )\  ,\
R^+(x)=\left (\begin{array}{cc}0&q^*(x)-\frac{\sqrt 2}{2}\\q(x)-\frac{\sqrt 2}{2}&0\end{array}\right )\  .$$
Then, using (\ref{matrix2.3}),  $\psi _1=X_1^++V_1$ where $V_1$ satisfies
$$i\partial _xV_1+M(Q^+-\lambda )V_1=-MR^+\psi _1\ ,$$
with $\ V_1(x,\lambda )\rightarrow 0$ as $x\rightarrow +\infty $,
hence, by the Duhamel formula, \beq\label{volt} \psi _1(x,\lambda
)=X_1^+(x,\lambda )+\int _x^{+\infty }S(x,y,\lambda )\psi
_1(y,\lambda )\, d\lambda \eeq where the matrix $S$ is given by \ben
\label{S} S(x,y,\lambda )=-ie^{i(x-y)M(Q^+-\lambda )}MR^+(y)\ .
 \een
Notice that, since $\lambda ^2\ge \frac12$ on $H$, the spectrum of
$M(Q^+-\lambda )$ is real. In view of the decay of $R^+(x)$ provided
by (\ref{estu}), $S$ satisfies
$$\forall a\in \R, \int _a^{+\infty }\sup _{x\geq a}|S(x,y,\lambda )|\, dy <+\infty .$$
Consequently, the  integral equation  ( \ref{volt}) is of Volterra
type, hence admits a unique solution $\psi _1$. We argue similarly
for the other Jost functions. Moreover, (\ref{estu}) also implies
that
$$\forall a\in \R, \int _a^{+\infty }\sup _{x\geq a}|\partial _\lambda ^jS(x,y,\lambda )|\, dy <+\infty ,\ j\leq 2.$$
Consequently, we can state
\begin{Lemma}\label{psiC2}
The Jost functions $\psi _1,\psi _2,\varphi _1,\varphi _2$ are $C^2$
functions on the hyperbola $H$ defined by (\ref{H}).
\end{Lemma}
\subsection{A representation formula}
In what follows, we are going to study further properties
of the Jost functions by establishing the following representation
formula, \ben\label{2.6}
\psi_1(x,\lambda)=X_1^+(x,\lambda)-\int_x^{+\infty}
\Psi(x,y)X_1^+(y,\lambda)dy, \een where $\Psi(x,y)$ is a
matrix-valued function. Substituting (\ref{2.6}) into
(\ref{matrix2.3}), we find that the matrix $\Psi(x,y)$ should
satisfy the following linear system
$$iM\partial _x\Psi +i\partial _y\Psi M-\Psi Q^++Q(x)\Psi =0$$
together with the boundary conditions
$$iM\Psi (x,x)-i\Psi (x,x)M+R^+(x)=0\ ,\  \int _x^{+\infty }|\Psi (x,y)|\, dy\rightarrow 0\ {\rm as}\ x\rightarrow +\infty .$$
A simple computation shows that this system is equivalent to
\beq\label{2.7}
\begin{split}
\Bigl(\f \p {\p x}+\f \p {\p y}\Bigr)\left(\begin{array}{c}
\Psi_{11}(x,y)\\
\Psi_{22}(x,y)
\end{array}\right)=i\left(\begin{array}{cc}
-\f{\sqrt{2}} 2 & q^*(x)\\
-q(x)& \f{\sqrt{2}} 2
\end{array}\right)\left(\begin{array}{c}
\Psi_{12}(x,y)\\
\Psi_{21}(x,y)
\end{array}\right),
\\
\Bigl(\f \p {\p x}-\f \p {\p y}\Bigr)\left(\begin{array}{c}
\Psi_{12}(x,y)\\
\Psi_{21}(x,y)
\end{array}\right)=i\left(\begin{array}{cc}
-\f{\sqrt{2}} 2 & q^*(x)\\
-q(x)& \f{\sqrt{2}} 2
\end{array}\right)\left(\begin{array}{c}
\Psi_{11}(x,y)\\
\Psi_{22}(x,y)
\end{array}\right),
\end{split}
\eeq
together with the boundary conditions
\begin{equation} \label{2.8}
\begin{split}
&\Psi_{12}^*(x,x)=\Psi_{21}(x,x)=-\f i 2(q(x)-\f {\sqrt{2}} 2),\\
& \int _x^{+\infty }|\Psi (x,y)|\, dy\rightarrow 0\ {\rm as}\ x\rightarrow +\infty.
 \end{split}
 \end{equation}
By the symmetry of the system (\ref{2.7}), we find that
\ben\label{symmetry} \Psi_{11}=\Psi_{22}^*,\qquad
\Psi_{12}=\Psi_{21}^*. \een From the invariance of the involution,
we have \beno X_2^+=\widetilde{X}_1^+, \qquad
\psi_2=\widetilde{\psi}_1. \eeno Hence, we can obtain a similar
representation for $\psi_2$ \ben\label{2.9}
\psi_2(x,\lambda)=X_2^+(x,\lambda)-\int_x^{+\infty}
\Psi(x,y)X_2^+(y,\lambda)dy. \een Similarly, we can also obtain the
following representation for $\varphi_1$ and $\varphi_2$
\ben\label{2.10}
\begin{split}
\varphi_1(x,\lambda)=X_1^-(x,\lambda)-\int_{-\infty}^x \Phi(x,y)X_1^-(y,\lambda)dy,\\
\varphi_2(x,\lambda)=X_2^-(x,\lambda)-\int_{-\infty}^x \Phi(x,y)X_2^-(y,\lambda)dy,
\end{split}
\een
where the matrix $\Phi(x,y)$ satisfies the linear system
\beq\label{2.11}
\begin{split}
\Bigl(\f \p {\p x}+\f \p {\p y}\Bigr)\left(\begin{array}{c}
\Phi_{11}(x,y)\\
\Phi_{22}(x,y)
\end{array}\right)=i\left(\begin{array}{cc}
\f{\sqrt{2}} 2 & q^*(x)\\
-q(x)& -\f{\sqrt{2}} 2
\end{array}\right)\left(\begin{array}{c}
\Phi_{12}(x,y)\\
\Phi_{21}(x,y)
\end{array}\right),
\\
\Bigl(\f \p {\p x}-\f \p {\p y}\Bigr)\left(\begin{array}{c}
\Phi_{12}(x,y)\\
\Phi_{21}(x,y)
\end{array}\right)=i\left(\begin{array}{cc}
\f{\sqrt{2}} 2 & q^*(x)\\
-q(x)& -\f{\sqrt{2}} 2
\end{array}\right)\left(\begin{array}{c}
\Phi_{11}(x,y)\\
\Phi_{22}(x,y)
\end{array}\right),
\end{split}
\eeq
together with the boundary conditions
\begin{equation} \label{2.12}
\begin{split}
&\Phi_{12}^*(x,x)=\Phi_{21}(x,x)=\f i 2(q(x)+\f {\sqrt{2}} 2),\\
&\int _{-\infty }^{x }| \Phi (x,y)|\, dy\rightarrow 0\ {\rm as}\ x\rightarrow -\infty\  .
 \end{split}
\end{equation}
Let us show how to solve the system (\ref{2.7}), (\ref{2.8}).
Introduce the nonnegative variable
$$p=\frac 12\, (y-x)\geq 0$$
and set
\beno V(x,p)=\left(\begin{array}{c}
\Psi_{12}(x,x+2p)\\
\Psi_{21}(x,x+2p)
\end{array}\right)\ , \ W(x,p)=\left(\begin{array}{c}
\Psi_{11}(x,x+2p)\\
\Psi_{22}(x,x+2p)
\end{array}\right)\ ,\ B(x)=i\left(\begin{array}{cc}
-\f{\sqrt{2}} 2 & q^*(x)\\
-q(x)& \f{\sqrt{2}} 2
\end{array}\right)
\eeno
Then system (\ref{2.7}) reads
$$\partial _xW(x,p)=B(x)V(x,p)\ ,\ \partial _pV(x,p)-\partial _xV(x,p)=-B(x)W(x,p)\ ,$$
and (\ref{2.8}) reads \beno V(x,0)=V_0(x)=\frac i2\left
(q(x)-\frac{\sqrt 2}2\right ) \left(\begin{array}{c}
1\\-1\end{array}\right)\ ,\ \int _0^{+\infty }(|V(x,p)|+|W(x,p) |)\,
dp \rightarrow 0\ {\rm as}\ x\rightarrow +\infty\ . \eeno Writing
$W$ from the first equation as \beq\label{Weq}
W(x,p)=-\int_x^{+\infty }B(x')V(x',p)\, dx'\ , \eeq we finally
obtain the integral equation for $V$, \beq \label{Veq}
V(x,p)=V_0(x+p)+\int _0^p\int
_{x+p-p'}^{+\infty}B(x+p-p')B(x')V(x',p')\, dx'\, dp'\ . \eeq Notice
that \ben\label{kernel} B(z)B(z')=-\left(\begin{array}{cc}
\f1 2-q(z')q^*(z) & \f{\sqrt{2}}2(q^*(z)-q^*(z'))\\
\f{\sqrt{2}}2(q(z)-q(z'))& \f1 2-q^*(z')q(z)
\end{array}\right)
\een
so that, in view of (\ref{estu}), if $z'\geq z\geq -a$ for a positive number $a$, then
$$|B(z)B(z')|\leq C(a)\langle z+a\rangle ^{-4}\ .$$
For every integer $N\geq 2$, denote by ${\mathcal V_N}$ the space of functions $V=V(x,p)$ on $\R \times \R _+$ such that,
for every $a>0$,
$$\sup _{x\geq -a}\langle x\rangle ^N\int _0^{+\infty }|V(x,p)|\, dp<+\infty $$
For every $V\in {\mathcal V_2}$, we set
$$Kv(x,p)=\int _0^p\int _{x+p-p'}^{+\infty }B(x+p-p')B(x')V(x',p')\, dx'\, dp'\ .$$
We claim that $Kv\in {\mathcal V_{N+1}}$ if $V\in {\mathcal V}_N$. Indeed, if $x\geq -a$,
$$|Kv(x,p)|\leq C(a)\int _0^p\langle x+p-p'+a\rangle ^{-4}\int _{x+p-p'}^{+\infty }|V(x',p')|\, dx'\, dp'\ ,$$
thus,
denoting  by $\| f\| $ the $L^1$ norm of $f=f(p)$ on $\R _+$,
we have
$$\| Kv(x,.)\|\leq C(a)\int _0^{+\infty }\langle x+s+a\rangle^{-4}\int _{x+s}^{+\infty}\| V(x',.)\| \, dx'\, ds$$
therefore
$$\| Kv(x,.)\| \leq C(a)\langle x+a\rangle ^{-3}\int _x^{+\infty }\| V(x',.)\| \, dx'.$$
In particular, if $\| V(x,.)\| \leq D(a)\langle x +a\rangle ^{-N}$, then
$$\| KV(x,.)\| \leq \frac{C(a)D(a)}{N-1}\langle x+a\rangle ^{-N-1}\ .$$
Starting from $V_0(x,p)=V_0(x+p)$, we conclude by an easy induction  that
$$\| K^nV_0(x,.)\| \leq \frac{C(a)^{n+1}}{(n+1)!}\langle x+a\rangle ^{-n-4}$$
which implies that the series $\sum _{n=0}^{+\infty }(-1)^nK^nV_0(x,.)$
converges in $L^1(\R _+)$ uniformly for $x\geq -a$, for all $a>0$.
Moreover, denoting by $V$ the sum of this series, $V$ solves (\ref{Veq}) and
$V\in {\mathcal V}_3$. Finally, coming back to (\ref{Weq}), we have $W\in {\mathcal V}_2$.
It is now a routine to estimate similarly the derivatives of order $k\leq 3$  of $V,W$ with respect
to $x,p$, and to show that they belong to ${\mathcal V}_3, {\mathcal V}_2$ respectively.
For future reference, we set
$$\underline \Psi (x,p)=\Psi (x,x+2p)\ ,$$
and we sum up the above results by \beq\label{Phi} \forall a>0,
\forall x\geq -a, \ \int _0^{+\infty}|\partial ^\alpha \underline
\Psi (x,p)|\, dp\leq C(a)\langle x+a\rangle ^{-3}\ , |\alpha |\leq
3. \eeq
\subsection{Analytic continuations and a priori bounds}
We are now in position to study the properties of Jost solutions
$\psi_1$ and $\psi_2$. Introduce the Riemann surface
$$\Gamma =\{ (\lambda ,\zeta )\in \C\times \C \ |\ \lambda ^2-\zeta
^2=\f 12\} \ ,$$ and denote by $\Gamma ^\pm $ the two sheets of
$\Gamma $ corresponding to $\pm {\rm Im}\zeta >0$. Notice that,
$\zeta $ is a single-valued holomorphic function of $\lambda $ on
$\Gamma ^+$ and on $\Gamma ^-$. Also notice that $$\overline \Gamma
^{\pm }=\Gamma ^{\pm }\cup H\ .$$

\begin{Lemma}\label{Lem2.3}
(1)\, The Jost solution $\psi_1(x,\lambda)\in C^3(\R_x)$ can be extended analytically to the lower sheet $\Gamma^-$
and has the form
\beno
\psi_1(x,\lambda)=X_1^+(x,\lambda )+\Psi_1(x,\zeta)X_1^+(x,\lambda).
\eeno

\no(2)\,\,The Jost solution $\psi_2(x,\lambda)\in C^3(\R_x)$ can be extended analytically to the upper sheet $\Gamma^+$
and has the form
\beno
\psi_2(x,\lambda)=X_2^+(x,\lambda )+\Psi_2(x,\zeta)X_2^+(x,\lambda).
\eeno
Here  $\Psi_1(x,\zeta)$ and $\Psi_2(x,\zeta)$ satisfy
\ben\label{2.30}
 |\Psi_1(x,\zeta)|+ |\Psi_2(x,\zeta)|
\le C(a)\langle \zeta\rangle ^{-1}\langle x+a\rangle ^{-3},\quad
\textrm{for}\,\, x\ge -a, a>0. \een
\end{Lemma}

\no Proof.\quad We rewrite (\ref{2.6}) as
\beno
\psi_1(x,\lambda)&=&X_1^+(x,\lambda)-\int_x^{+\infty} \Psi(x,y)X_1^+(y,\lambda)dy
\\
&=&X_1^+(x,\lambda)-2\int_0^{+\infty}\underline \Psi(x,p)e^{-2i\zeta p}\, dp\,  X_1^+(x,\lambda)\\
&\triangleq&  X_1^+(x,\lambda)+ \Psi_1(x,\zeta)X_1^+(x,\lambda).
\eeno
We get by integrating by parts that
\beno
\Psi_1(x,\zeta)=\f i \zeta\underline \Psi(x,0)
+\f i \zeta\int_0^{+\infty}\p_p\underline \Psi(x,p)e^{-2i\zeta p}dp\  ,
\eeno
which, in view of the estimates (\ref{Phi}) on $\underline \Psi $, implies (\ref{2.30}).
Since $(2)$ can be proved in a similar way, we omit its proof.\endproof

Similarly, we state the following properties for the Jost solutions $\varphi_1$ and $\varphi_2$.

\begin{Lemma}\label{Lem2.4}
(1)\, The Jost solution $\varphi_1(x,\lambda)\in C^3(\R_x)$ can be extended analytically to the upper sheet $\Gamma^+$
and has the form
\beno
\varphi_1(x,\lambda)=X_2^{-}(x,\lambda)+{\Phi}_1(x,\zeta)X_1^-(x,\lambda).
\eeno

\no(2)\,\,The Jost solution $\varphi_2(x,\lambda)\in C^3(\R_x)$ can
be extended analytically to the lower sheet $\Gamma^-$ and has the
form \beno \varphi_2(x,\lambda)=X_2^{-}(x,\lambda
)+{\Phi}_2(x,\zeta)X_2^-(x,\lambda). \eeno Here ${\Phi}_1(x,\zeta)$
and ${\Phi}_2(x,\zeta)$ satisfy \ben\label{2.31}
|{\Phi}_1(x,\zeta)|+ |{\Phi}_2(x,\zeta)| \le C(a)\langle
\zeta\rangle ^{-1}\langle x-a\rangle ^{-3},\quad \textrm{for}\,\,
x\leq a, a>0. \een
\end{Lemma}

\subsection{The unperturbed case}
As a next step we  apply the above results to the particular case
$$q^0(x)=\frac{\sqrt 2}{2}\tanh \left (\frac{x}{\sqrt 2}\right )\ .$$
In this  case, since $q^0$ is real valued, the kernel of equation
(\ref{Veq}) satisfies, in view of (\ref{kernel}),
\beno
B(z)B(z')\left (\begin{array}{c}1\\-1\end {array}\right )&=-\left (\frac 12-q(z)q(z')-\frac {\sqrt 2}{2}q(z)+\frac {\sqrt 2}{2}q(z')\right )\left (\begin{array}{c}1\\-1\end {array}\right )\\
&=-\frac 12\left (1-\tanh \left (\frac {z}{\sqrt 2}\right )\right )
\left (1+\tanh\left (\frac{z'}{\sqrt 2}\right )\right )\left (\begin{array}{c}1\\-1\end {array}\right )\ ,
\eeno
so that $V_2(x,p)=-V_1(x,p)$ and (\ref{Veq}) reads
\beno\label{V1eq}
V_1(x,p)=-\frac{i}{2\sqrt 2}\left (1-\tanh \left (\frac{x+p}{\sqrt 2}\right )\right )\hskip 8cm\\-\frac 12\int _0^p\int _{x+p-p'}^{+\infty }
\left (1-\tanh \left (\frac {x+p-p'}{\sqrt 2}\right )\right )
\left (1+\tanh\left (\frac{x'}{\sqrt 2}\right )\right )\, V_1(x',p')\, dx'\, dp'\ .
\eeno
Introducing the following new variables
$$s=e^{\sqrt 2(x+p)}>0 \ ,\ r=e^{\sqrt 2p}\geq 1\ ,\ v(s,r)=V_1(x,p)\ ,$$
we obtain
$$v(s,r)=-\frac{i}{2\sqrt 2(s+1)}-\int _1^r\int _s^{+\infty }\frac{v(s',r')}{(s+r')(s'+r')}\, ds'\, dr'$$
and it is easy to see that this integral equation admits for solution
$$v(s,r)=-\frac{i}{2\sqrt 2(s+r)}\ .$$
Coming back to the original variables, we infer \ben\label{Psi0}
\begin{split}\label{Psi0}
\Psi _{12}^0(x,x+2p)=V_1^0(x,p)=-\frac{ie^{-\sqrt 2p}}{\sqrt
2(1+e^{\sqrt
2x})}\ ,\\
\Psi _{11}^0(x,x+2p)=W_1^0(x,p)=\frac{e^{-\sqrt
2p}}{\sqrt 2(1+e^{\sqrt 2x})}\ .
\end{split}
\een
 Through the representation derived in Section 2.1, we
obtain \ben\label{2.15}
\begin{split}
\psi_1^0(x,\lambda)=e^{-i\zeta x}\left(\begin{array}{c}
1-\f 1 {1+e^{\sqrt{2}x}}\f {\f {\sqrt{2}}2-i(\lambda-\zeta)} {\f {\sqrt{2}}2+i\zeta}\\
\sqrt{2}(\lambda-\zeta)-\f 1 {1+e^{\sqrt{2}x}}\f {\f {\sqrt{2}}2i+(\lambda-\zeta)} {\f {\sqrt{2}}2+i\zeta}
\end{array}\right),\\
\psi_2^0(x,\lambda)=e^{i\zeta x}\left(\begin{array}{c}
\sqrt{2}(\lambda-\zeta)-\f 1 {1+e^{\sqrt{2}x}}\f {-\f {\sqrt{2}}2i+(\lambda-\zeta)} {\f {\sqrt{2}}2-i\zeta}\\
1-\f 1 {1+e^{\sqrt{2}x}}\f {\f {\sqrt{2}}2+i(\lambda-\zeta)} {\f {\sqrt{2}}2-i\zeta}
\end{array}\right),
\end{split}
\een
and similarly, through Lemma \ref{Lem2.4},
\ben\label{2.16}
\begin{split}
&\varphi_1^0(x,\lambda)=e^{-i\zeta x}\left(\begin{array}{c}
1-\f {e^{\sqrt{2}x}} {1+e^{\sqrt{2}x}}\f {\f {\sqrt{2}}2+i(\lambda-\zeta)} {\f {\sqrt{2}}2-i\zeta}\\
-\sqrt{2}(\lambda-\zeta)+\f {e^{\sqrt{2}x}} {1+e^{\sqrt{2}x}}\f {-\f {\sqrt{2}}2i+(\lambda-\zeta)} {\f {\sqrt{2}}2-i\zeta}
\end{array}\right),\\
&\varphi_2^0(x,\lambda)=e^{i\zeta x}\left(\begin{array}{c}
-\sqrt{2}(\lambda-\zeta)+\f {e^{\sqrt{2}x}} {1+e^{\sqrt{2}x}}\f {\f {\sqrt{2}}2i+(\lambda-\zeta)} {\f {\sqrt{2}}2+i\zeta}\\
1-\f {e^{\sqrt{2}x}} {1+e^{\sqrt{2}x}}\f {\f {\sqrt{2}}2-i(\lambda-\zeta)} {\f {\sqrt{2}}2+i\zeta}
\end{array}\right).
\end{split}
\een

\subsection{Perturbation analysis}
We close this section by describing how the  results of the previous
subsections can be precised when $q$ is assumed to be a perturbation
of $q^0$ in the sense of Theorem \ref{1.1}, namely if
\begin{equation}\label{qperturb}
q(x)=q^0(x)+\varepsilon q^1(x)\ ,\ \langle x\rangle ^4\left
|\frac{d^kq^1}{dx^k}(x)\right |\leq 1\ ,\ 0\leq k\leq 3\ ,
\end{equation}
and $\varepsilon $ is a small positive parameter. Revisiting the
analysis of subsections 2.1 and 2.2, and comparing to the results of
subsection 2.3, the following perturbation results can be easily
proved.
\begin{Lemma}\label{PsiPhiperturb}
The matrix-valued  functions $\Psi $, $\Phi $ defined by
(\ref{2.7}), (\ref{2.8}), (\ref{2.11}), (\ref{2.12}) can be written
as
$$\Psi (x,y)=\Psi ^0(x,y)+\vep \Psi ^1(x,y)\ ,\ \Phi (x,y)=\Phi ^0(x,y)+\vep \Phi
^1(x,y)\ ,$$ with the following estimates : if we set, for $p\ge 0$,
$$\underline \Psi ^1 (x,p):=\Psi ^1(x,x+2p)\ ,\ \underline \Phi ^1
(x,p)=\Phi ^1(x,x-2p)\ ,$$ then \ben\label{estPsiPhiperturb} \forall
a>0, \forall x\geq -a, \ \int _0^{+\infty}|\partial ^\alpha
\underline \Psi ^1(x,p)|\, dp\leq C(a)\langle x+a\rangle ^{-3}\ ,
|\alpha |\leq 3\ ,\\
\forall a>0, \forall x\leq a, \ \int _0^{+\infty}|\partial ^\alpha
\underline \Phi ^1(x,p)|\, dp\leq C(a)\langle x-a\rangle ^{-3}\ ,
|\alpha |\leq 3\ .\een
\end{Lemma}

\begin{Lemma}\label{psiphiperturb}
For fixed $q^1$, the Jost functions are real analytic of the
paramter $\varepsilon $ in a neighborhood of $0$ in $\R $. Moreover,
for $\varepsilon$ small enough, they  can be written as
 \beno
\psi_1(x,\lambda)=\psi_1^0(x,\lambda)+\varepsilon
\psi_1^1(x,\lambda)\ ,\
\psi_2(x,\lambda)=\psi_2^0(x,\lambda)+\varepsilon
\psi_2^1(x,\lambda)\\
\varphi_1(x,\lambda)=\varphi_1^0(x,\lambda)+\vep
{\varphi}_1^1(x,\lambda)\ ,\
\varphi_2(x,\lambda)=\varphi_2^0(x,\lambda)+\vep
{\varphi}_2^1(x,\lambda), \eeno with \beno
 |\psi_1^1(x,\lambda)|+|\psi_2^1(x,\lambda)|
\le C\langle \zeta\rangle ^{-1}\langle x\rangle ^{-3},\quad \textrm{for}\,\, x\ge 0\\
|{\varphi}_1^1(x,\lambda)|+ |{\varphi}_2^1(x,\lambda)| \le
C\langle\zeta\rangle ^{-1}\langle x\rangle ^{-3},\quad
\textrm{for}\,\, x\le 0,\eeno and $\lambda \in \Gamma ^-$ in the
cases of $\psi_1,\varphi_2$, while $\lambda \in \Gamma ^+$ in the
cases of $\psi_2,\varphi_1$.
\end{Lemma}

\section{Transition coefficients and their properties}

Now let $\zeta$ be real. Let $\psi_1(x,\lambda), \psi_2(x,\lambda),
\varphi_1(x,\lambda),$ and $\varphi_2(x,\lambda)$ be the Jost
solutions constructed in Section 2. Since their Wronskian is
independent of $x$, we get by their asymptotic behaviour at infinity
that \ben
\{\psi_1,\psi_2\}=\{\varphi_1,\varphi_2\}=4\zeta(\lambda-\zeta).\label{3.1}
\een Thus, $\psi_1$ and $\psi_2$ are linearly independent for
$\zeta\neq 0$. Hence, we can expand $\varphi_1$ and $\varphi_2$
uniquely as \ben \varphi_1=a\psi_1+b\psi_2,\qquad
\varphi_2=a^*\psi_2+b^*\psi_1,\label{3.2} \een from which and
(\ref{3.1}), we get \ben a(\lambda)=\f {\{\varphi_1, \psi_2\}}
{4\zeta(\lambda-\zeta)},\quad b(\lambda)=-\f {\{\varphi_1, \psi_1\}}
{4\zeta(\lambda-\zeta)},\label{3.3} \een which are called the
transition coefficients. In the case when $q(x)=\f {\sqrt{2}}
2\tanh(\f x {\sqrt{2}})$, it is easy to compute by (\ref{2.15}) and
(\ref{2.16}) that the corresponding transition coefficients denoted
by $a^0(\lambda)$ and $b^0(\lambda)$ are
$$
a^0(\lambda)=\f {\lambda+\zeta-\f {\sqrt{2}} 2i} {\lambda+\zeta+\f
{\sqrt{2}} 2i},\quad b^0(\lambda)=0.
$$
Note that
$$
\{\varphi_1,\varphi_2\}=(|a|^2-|b|^2)\{\psi_1,\psi_2\},
$$
and $\{\varphi_1,\varphi_2\}=\{\psi_1,\psi_2\}$, we have the
normalization relation \ben |a(\lambda)|^2-|b(\lambda)|^2=1, \quad
\textrm{for}\,\, \lambda\in \R, \,|\lambda|>\f {\sqrt{2}} 2.
\label{3.4} \een We find by Lemma \ref{Lem2.3}-\ref{Lem2.4} and the
definition of $a(\lambda)$ that the function $a(\lambda)$ can be
extended analytically to the upper sheet $\Gamma^+$ . In what
follows, we will study the properties about the transition
coefficients. Let us begin with the following simple fact about the
function $\zeta(\lambda)$ which will be constantly used: For
$\lambda\in \R, |\lambda|\gg 1$ \ben\label{3.5} \zeta=\lambda+O(\f 1
{\lambda}), \,\,\textrm{for}\,\, \lambda>0, \qquad
\zeta=-\lambda+O(\f 1 {\lambda}),\,\, \textrm{for}\,\, \lambda<0.
\een

\begin{Lemma}\label{Lem3.1} Let $a(\lambda), b(\lambda)$ be given by (\ref{3.3}).
Then there hold for $\lambda\in \R, |\lambda|>\f {\sqrt{2}} 2$
\ben
&&|a(\lambda)-a^0(\lambda)|\le C\vep|\zeta|^{-1},\label{3.6}\\
&&|b(\lambda)|\le C\vep\min(|\zeta|^{-1},|\zeta|^{-3}),\label{3.7}\\
&&\left |\f {b(\lambda)} {a(\lambda)}\right |\le
\min(1,C\vep|\zeta|^{-1},C\vep|\zeta|^{-3}).\label{3.8} \een
\end{Lemma}

\no Proof.\quad We write \ben \psi_2=\psi_2^0+\vep \psi_2^1,\quad
\varphi_1=\varphi_1^0+\vep \varphi_1^1\ .\label{3.9} \een
Substituting (\ref{3.9}) into the first formula of (\ref{3.3}), we
get \beno a(\lambda)&=&\f {\{\varphi_1,\psi_2\}}
{4\zeta(\lambda-\zeta)}=
\f {\{\varphi_1^0+\vep \varphi_1^1,\psi^0_2+\vep \psi^1_2\}} {4\zeta(\lambda-\zeta)}\nonumber\\
&=& \f {\{\varphi_1^0,\psi^0_2\}} {4\zeta(\lambda-\zeta)}+
\vep \f {\{\varphi_1^0,\psi^1_2\}+\{ \varphi_1^1,\psi^0_2\}+\vep\{ \varphi_1^1, \psi^1_2\}} {4\zeta(\lambda-\zeta)}\nonumber\\
&\triangleq& a^0(\lambda)+\vep a^1(\lambda), \eeno which together
with Lemma \ref{psiphiperturb} and (\ref{3.5}) gives (\ref{3.6}). We
next prove (\ref{3.7}). Note that $\{\varphi_1^0, \psi_1^0\}=0$, we
get \ben\label{3.10} \{\varphi_1,\psi_1\}=\vep
\{\varphi_1^0,\psi^1_1\}+\vep \{\varphi_1^1,\psi^0_1\}+\vep
^2\{\varphi_1^1,\psi^1_1\}. \een Thus, we can get by  using Lemma
\ref{psiphiperturb} and (\ref{3.5}) that \ben\label{3.11}
|b(\lambda)|\le C\vep|\zeta|^{-1}. \een In order to obtain the
better decay estimate for $b(\lambda)$, we need to use the following
more subtle argument. Multiplying by $e^{2i\zeta x}$ on both sides
of (\ref{3.10}) and taking the derivative with respect to $x$ to the
resulting equation, we get by using the fact that
$\{\varphi_1,\psi_1\}$ is independent of $x$ that \ben\label{3.12}
-8i\zeta^3e^{2i\zeta x}\{\varphi_1,\psi_1\}=\vep\bigl(\f {d}
{dx}\bigr)^3\bigl(\{e^{i\zeta x}\varphi_1^0,e^{i\zeta
x}\psi^1_1\}+\{e^{i\zeta x}\varphi_1^1,e^{i\zeta x}\psi^0_1\} +\nonumber \\
+\vep \{e^{i\zeta x}\varphi_1^1,e^{i\zeta x}\psi^1_1\}\bigr). \een
It is easy to verify by (\ref{2.16}) that for $k\in \N$
\ben\label{3.13} |\bigl(\f {d} {dx}\bigr)^k(e^{i\zeta
x}\varphi_1^0(x,\lambda))|\le C. \een By the proof of Lemma
\ref{Lem2.3} adapted to Lemma \ref{psiphiperturb}, we have
\ben\label{3.14} |\bigl(\f {d} {dx}\bigr)^k(e^{i\zeta
x}\psi^1_1(x,\zeta))|\le C,\quad \textrm{for} \quad k\le 3, x\ge 0.
\een By summing up (\ref{3.13}) and (\ref{3.14}), we obtain \beno
|\bigl(\f {d} {dx}\bigr)^3\{e^{i\zeta x}\varphi_1^0,e^{i\zeta
x}\psi^1_1\}|\le C\ . \eeno Similarly, it can be proved that the
other two terms on the right side of (\ref{3.12}) are bounded by $C$
. Thus, we get by using (\ref{3.5}) and (\ref{3.12}) that
$$
|b(\lambda)|\le C\varepsilon|\zeta|^{-3},\qquad for \quad |\zeta|\gtrsim 1,
$$
which together with (\ref{3.11}) gives (\ref{3.7}). Now we turn to
prove (\ref{3.8}). We get by (\ref{3.4}) that \beno \left |\f
{b(\lambda)} {a(\lambda)}\right | \le 1,\qquad |a(\lambda)|\ge 1,
\eeno from which and (\ref{3.7}), it follows (\ref{3.8}).\endproof

\begin{Lemma}\label{Lem3.2}
The function $a(\lambda)$ has at most one zero. Assume that $\lambda_0$ is a zero of $a(\lambda)$, one has

\noindent (1)\, $\lambda_0\in (-\f {\sqrt{2}} 2, +\f {\sqrt{2}} 2)$ is simple;

\noindent (2)\, there exists a constant $b_0$ such that $\varphi_1(x,\lambda_0)=b_0\psi_2(x,\lambda_0)$ for any $x\in \R$;

\noindent (3)\, for $\varepsilon$ small enough, there holds
\ben\label{3.15}
|\lambda_0|\le C\varepsilon,\quad
|b_0-i|\le C\varepsilon,\quad
|\mu_0+2|\le C\vep,
\een
where
$$
\mu_0=\f {b_0} {\nu_0a'(\lambda_0)}, \quad \nu_0=(\f12-\lambda_0^2)^\f12.
$$
\end{Lemma}

\no Proof.\,\,Let us assume that $\lambda_0$ is a zero of
$a(\lambda)$. From (\ref{3.4}), we see that $\lambda $ does not
belong to $\{ \lambda\in \R, |\lambda|>\f {\sqrt{2}} 2\} $, hence
$Im\zeta >0$.

Since
$$
\{\varphi_1,\psi_2\}(\lambda_0)=0,
$$
there exists a constant $b_0$ such that \beno
\varphi_1(x,\lambda_0)=b_0\psi_2(x,\lambda_0),\quad x\in \R, \eeno
which implies that the system (\ref{2.3}) with $\lambda=\lambda_0$
has a global $L^2$ solution  on $\R $. Thus, the zeros of
$a(\lambda)$ correspond to the eigenvalues of the system
(\ref{2.3}). From the self-adjoint character of the system , it
follows that $\lambda_0$ must lie on the segment $(-\f {\sqrt{2}} 2,
+\f {\sqrt{2}} 2)$ of the real axis. We next prove that $\lambda_0$
is simple. It suffices to prove that $a'(\lambda_0)\neq 0.$ By
(\ref{3.3}), we have \beno a'(\lambda_0)=\f {\{\f {\partial}
{\partial \lambda}\varphi_1,\psi_2\}(x,\lambda_0) +\{\varphi_1,\f
{\partial} {\partial \lambda}\psi_2\}(x,\lambda_0)}
{4\zeta_0(\lambda_0-\zeta_0)},\quad \zeta_0=(\lambda_0^2-\f12)^\f12
=i\nu _0\ . \eeno Set $ \sigma=\left(\begin{array}{cc} 1&0\\0&-1
\end{array}\right),
$
we get by using the equations (\ref{2.3}) that
\beno
\f {\partial} {\partial x}\{\f {\partial} {\partial \lambda}\varphi_1,\psi_2\}(x,\lambda)
&=&-i\{\sigma \varphi_1, \psi_2\}(x,\lambda),\\
\f {\partial} {\partial x}\{\varphi_1,\f {\partial} {\partial \lambda}\psi_2\}(x,\lambda)
&=&i\{\sigma \varphi_1, \psi_2\}(x,\lambda),
\eeno
from which and the fact that the Jost solutions $\varphi_1(x,\lambda_0)$ and $\psi_2(x,\lambda_0)$
decay to zero as $|x|\rightarrow \infty$, it follows that
\beno
\{\f {\partial} {\partial \lambda}\varphi_1,\psi_2\}(x,\lambda_0)
=-ib_0\int_{-\infty}^x\{\sigma \psi_2, \psi_2\}(x',\lambda_0)dx',\\
\{\varphi_1,\f {\partial} {\partial \lambda}\psi_2\}(x,\lambda_0)
=-ib_0\int^{+\infty}_x\{\sigma \psi_2, \psi_2\}(x',\lambda_0)dx'.
\eeno
From the invariance of the involution of (\ref{2.3}) and the uniqueness of the Jost solutions , we find that
for $\lambda\in (-\f {\sqrt{2}} 2, +\f {\sqrt{2}} 2)$
$$
\sqrt{2}(\lambda+\zeta)\psi_2(x,\lambda)=\widetilde{\psi}_2(x,\lambda),
$$
which implies that
$$
\{\sigma \psi_2, \psi_2\}(x,\lambda)=|\psi_2(x,\lambda)|^2.
$$
Putting together the above formulas, we get
\ben\label{3.16}
a'(\lambda_0)=-\f {ib_0\int_{-\infty}^{+\infty}|\psi_2(x,\lambda_0)|^2dx} {2\sqrt{2}\zeta_0}\neq 0.
\een

We now prove (\ref{3.15}). Using the similar proof as in (\ref{3.6}), we have
$$
|a(\lambda)-a^0(\lambda)|\le C\varepsilon|\zeta|^{-1},\quad
\textrm{for}\,\,\lambda \in (-\f {\sqrt{2}} 2, +\f {\sqrt{2}} 2).
$$
Set $\zeta=i(\f12 -\lambda^2)^\f12=i\nu$ for $\lambda \in (-\f {\sqrt{2}} 2, +\f {\sqrt{2}} 2)$. Then we have
$$
|a^0(\lambda)|=\f {\sqrt{2}|\lambda|} {1+\sqrt{2}\nu}\ge \f
{\sqrt{2}} 2|\lambda|.
$$
Hence, we have at the zero point $\lambda_0$ that
\beno
|\lambda_0|^2|\zeta_0|^2\le C\varepsilon^2,\quad \zeta_0=(\lambda^2-\f12)^\f12
\eeno
which is equivalent to the inequality
$$
\lambda_0^4-\f12 \lambda_0^2+C\varepsilon^2\ge 0.
$$
We get by solving the above inequality that
\ben\label{3.17}
\lambda_0^2\in \Bigl[0,\f14-\f12\sqrt{\f14-4C\varepsilon^2}\Bigr]\cup\Bigl[\f14+\f12\sqrt{\f14-4C\varepsilon^2},\f12\Bigr).
\een
Note that $\f14+\f12\sqrt{\f14-4C\varepsilon^2}\rightarrow \f12$ as $\varepsilon\rightarrow 0$, while by (\ref{3.4})
$$
\liminf_{\lambda\rightarrow \pm\f12} |a(\lambda)|\ge 1,
$$
which together with (\ref{3.17}) implies that for $\varepsilon$ small enough
\ben\label{3.18}
|\lambda_0|\le c\vep.
\een

Using the fact that $\varphi_1^0(x,0)=i\psi_2^0(x,0)$, we get \beno
(i-b_0)\psi_2^0(x,0)&=&\varepsilon
(b_0\psi_2^1(x,\lambda_0)-\varphi_1^1(x,\lambda_0))
-(\varphi_1^0(x,\lambda_0)-\varphi_1^0(x,0))\\
&&+b_0(\psi_2^0(x,\lambda_0)-\psi_2^0(x,0)),
\eeno
which together with Lemma \ref{psiphiperturb} and (\ref{3.18}) gives
$$
|b_0-i||\psi_2^0(0,0)|=\f {\sqrt{2}}2|b_0-i|\le C\varepsilon|b_0-i|+C\varepsilon,
$$
which implies the second inequality of (\ref{3.15}). We next prove the third inequality of (\ref{3.15}).
By the first two inequalities of (\ref{3.15}) and (\ref{3.16}), it suffices to prove that
\ben
\bigl|\int_{-\infty}^{+\infty}|\psi_2(x,\lambda_0)|^2dx-\sqrt{2}\bigr|\le C\varepsilon.\label{3.19}
\een
We write
\beno
\int_{-\infty}^{+\infty}|\psi_2(x,\lambda_0)|^2dx
&=&\int_{0}^{+\infty}|\psi_2^0(x,\lambda_0)+\varepsilon\psi_2^1(x,\lambda_0)|^2dx\\
&&+\f 1 {|b_0|^2}\int_{-\infty}^{0}|\varphi_1^0(x,\lambda_0)+\varepsilon\varphi_1^1(x,\lambda_0)|^2dx,
\eeno
from which and Lemma \ref{psiphiperturb}, it follows that
\ben\label{3.20}
\bigl|\int_{-\infty}^{+\infty}|\psi_2(x,\lambda_0)|^2dx-\int_{0}^{+\infty}|\psi_2^0(x,\lambda_0)|^2dx
-\f 1 {|b_0|^2}\int_{-\infty}^{0}|\varphi_1^0(x,\lambda_0)|^2dx\bigr|\le C\varepsilon.
\een
On the other hand, from the exact formula (\ref{2.15}-\ref{2.16}) and (\ref{3.18}), it is easy to verify that
\beno
&&\varphi_1^0(x,0)=i\psi_2^0(x,0),\quad \int_{-\infty}^{+\infty}|\psi_2^0(x,0)|^2dx=\sqrt{2}, \quad \textrm{and}\\
&& \int_{0}^{+\infty}|\psi_2^0(x,\lambda_0)-\psi_2^0(x,0)|^2dx+
\int_{-\infty}^{0}|\varphi_1^0(x,\lambda_0)-\varphi_1^0(x,0)|^2dx\le C\varepsilon,
\eeno
which together with (\ref{3.20}) imply (\ref{3.19}).

Finally, we prove that $a(\lambda)$ has at most one zero. Assume that $\lambda_0'$ is another zero of $a(\lambda)$.
Then $|\lambda_0'|\le C\vep$ by the above proof.
We get by an exact computation that
$$
a^0(\lambda_0)-a^0(\lambda_0')=\f {i\sqrt{2}}
{(1+\sqrt{2}\nu_0)(1+\sqrt{2}\nu_0')}
\Bigl((1+\sqrt{2}\nu_0)(\lambda_0'-\lambda_0)+\sqrt{2}\lambda_0(\nu_0-\nu_0')\Bigr).
$$
Hence, for $\varepsilon$ small enough there exists a constant $c>0$
such that \ben |a^0(\lambda_0)-a^0(\lambda_0')|\ge
c|\lambda_0-\lambda_0'|.\label{3.21} \een Using Lemma
\ref{psiphiperturb} and the fact that
$|\lambda_0|,|\lambda_0'|\le C\varepsilon$, we have \beno
|a^1(\lambda_0)-a^1(\lambda_0')|\le C|\lambda_0-\lambda_0'|, \eeno
which together with (\ref{3.21}) implies that
$\lambda_0=\lambda_0'$. The proof of Lemma \ref{Lem3.2} is
complete.\endproof

\section{The evolution of the transition coefficients}

In this section, we will derive the evolution equations of the transition coefficients.
Let us begin by deriving the evolution equations of the Jost solutions.

Differentiating (\ref{2.1}) with respect to $t$, we get
$$
\f \p {\p t}L_u\chi+L_u\f {\p \chi} {\p t}=E\chi,
$$
from which and (\ref{4.2}), it follows that
$$
L_u(\f{\p \chi} {\p t}+iB_u\chi)=E(\f{\p \chi} {\p t}+iB_u\chi),
$$
which together with the asymptotic behavior and uniqueness of the Jost solution gives
\beq\label{4.3}
\f{\p \chi} {\p t}+iB_u\chi=i\sqrt{3}(\f E 2-\zeta)^2\chi, \quad
\eeq
for $\chi$ determined by the Jost solution $\psi_1$ or $\varphi_1$, and
\beq\label{4.4}
\f{\p \chi} {\p t}+iB_u\chi=i\sqrt{3}(\f E 2+\zeta)^2\chi, \quad
\eeq
for $\chi$ determined by the Jost solution  $\psi_2$ or $\varphi_2$.

Differentiating (\ref{3.2}) with respect to $t$, we get by using (\ref{4.3}) and (\ref{4.4}) that
\beno
\f {d a} {dt}\psi_1+(\f {d b} {dt}+i4\lambda\zeta b)\psi_2=0.
\eeno
So, we get the evolution equations of the transition coefficients
\ben\label{4.5}
a(t,\lambda)=a(0,\lambda),\qquad b(t,\lambda)=b(0,\lambda)\exp(-i4\lambda\zeta t).
\een
Thus, the zero of $a(t,\lambda)$ does not depend on the time $t$.

At the zero $\lambda_0$ of $a(\lambda)$, we have
$$
\varphi_1(x,\lambda_0)=b_0\psi_2(x,\lambda_0).
$$
We get by using the evolution equations (\ref{4.3}) and (\ref{4.4}) of the Jost solutions  that
$$
\f {d b_0(t)} {dt}+i4\lambda_0\zeta_0 b_0(t)=0,
$$
which gives
\ben\label{4.6}
b_0(t)=b_0(0)\exp(-i4\lambda_0\zeta_0 t), \quad \zeta_0=(\lambda_0^2-\f12)^\f12.
\een

\section{The Marchenko equations}

In this section, we will derive the Marchenko equations which
connect the potential $q(x)$ with the scattering data. Throughout
this section, we assume that the function $a(\lambda)$ has one zero
$\lambda_0$, and we denote by  $b_0$ the constant such that \ben
\varphi_1(x,\lambda_0)=b_0\psi_2(x,\lambda_0),\label{5.1} \een and
$$\nu _0=\sqrt {\f 12-\lambda _0^2}\ ,\ \mu _0=\f {b_0} {\nu_0a'(\lambda_0)}.$$
From (\ref{3.16}), we know that $\mu_0$ is a real number. Moreover,
by the the uniqueness of the Jost solution, the change $(\lambda
,\zeta )\in H\mapsto (\lambda ,-\zeta )\in H$ implies that
$$
\varphi_1\mapsto \f {1} {\sqrt{2}(\lambda-\zeta)}\varphi_2,\qquad
\psi_1\mapsto \f 1 {\sqrt{2}(\lambda-\zeta)}\psi_2,
$$
which together with (\ref{3.2}) implies that
\ben\label{5.1a}
a\mapsto -a^*,\qquad b\mapsto -b^*.
\een
For $\zeta \in \R$, we set
$$\lambda =\lambda (\zeta )=\sqrt {\zeta ^2 +\f 12}\ , $$
and
$$c_1(\zeta)= c(\lambda)+c(-\lambda),\qquad c_2(\zeta)=\f {c(\lambda)-c(-\lambda)}{\lambda},\quad
c=\f b a.$$ The above observations imply that
$$c_1(-\zeta )=c_1(\zeta )^*\ ,\ c_2(-\zeta )=c_2(\zeta )^* \ .$$
We then define the following real valued functions of the real
variable $z$, \ben \label{F} F_1(z)=\f
1{2\pi}\int_{-\infty}^{+\infty}c_1(\zeta)e^{i\zeta z}d\zeta-\mu
_0\lambda _0e^{-\nu_0z}\ ,\  F_2(z)=\f
1{2\pi}\int_{-\infty}^{+\infty}c_2(\zeta)e^{i\zeta z}d\zeta -\mu_0
e^{-\nu_0z}\ . \een
\begin{Proposition}\label{marchenko}
The functions $\Psi _{11}$, $\Psi _{12}$ defined by (\ref{2.7}),
(\ref{2.8}) satisfy the following system for $y\ge x$,
\ben\label{5.7}
&&2\sqrt{2}\Psi_{11}(x,y)=F_2(x+y)\nonumber\\
&&\quad-\int_x^{+\infty}\Psi_{12}(x,s)\sqrt{2}(F_1(s+y)-iF_2'(s+y))+\Psi_{11}(x,s)F_2(s+y)ds\
,
\een
and
\ben\label{5.8}
&& 2\sqrt{2}\Psi_{12}(x,y)=\sqrt{2}(F_1(x+y)+iF_2'(x+y))\nonumber\\
&&\quad-\int_x^{+\infty}\Psi_{11}(x,s)\sqrt{2}(F_1(s+y)+iF_2'(s+y))+\Psi_{12}(x,s)F_2(s+y)ds\
. \een
\end{Proposition}
Let us prove Proposition \ref{marchenko}. We rewrite the first
equation of (\ref{3.2}) in the form \beq\label{5.2} \left (\f
{\varphi_1} a-X_1^+\right )\f {e^{i\zeta y}}{2\pi \zeta}=\left
(\psi_1-X_1^++\f ba\psi_2\right )\f{e^{i\zeta y}}{2\pi\zeta},\quad
y\ge x. \eeq The left-hand side of (\ref{5.2}) is analytic on the
upper sheet $\Gamma^+$ of the Riemann surface $\Gamma$, with the
exception of the point $\lambda_0$. We integrate (\ref{5.2}) along
the contour with respect to $\lambda$ indicated in the following
figure:
\begin{figure}[H]
 \centering
\scalebox{0.8}{\includegraphics{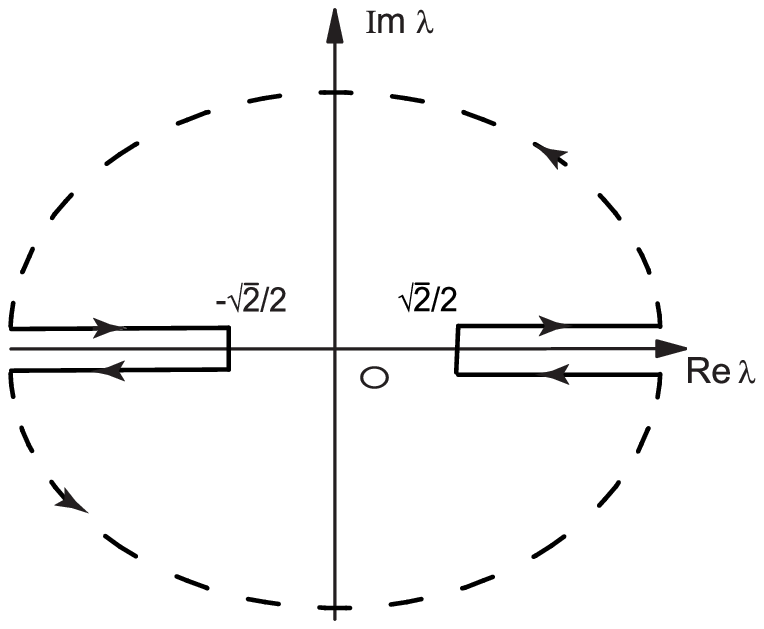}}
\end{figure}
\no We get by using the residue theorem and (\ref{5.2}) that
\ben\label{5.3} \int \left (\f1a\varphi_1-X_1^+\right )e^{i\zeta
y}\f {d\lambda} {2\pi \zeta }=\f {\varphi _1(x,\lambda_0)}
{\nu_0a'(\lambda_0)} e^{-\nu_0 y} =b_0\f {\psi_2(x,\lambda_0)}
{\nu_0a'(\lambda_0)}e^{-\nu_0 y}. \een Substituting the
representation (\ref{2.9}) into (\ref{5.3}), we get \ben\label{5.4}
&&\int\left (\f1a\varphi_1-X_1^+\right )e^{i\zeta
y}\f{d\lambda}{2\pi \zeta}=\left(\begin{array}{c}
\sqrt{2}F_1^{(1)}(x+y)+i\sqrt{2}F_2^{(1)'}(x+y)\\
F_2^{(1)}(x+y)
\end{array}\right)\nonumber\\
&&\qquad\qquad-\int_x^\infty\Psi(x,s)\left(\begin{array}{c}
\sqrt{2}F_1^{(1)}(s+y)+i\sqrt{2}F_2^{(1)'}(s+y)\\
F_2^{(1)}(s+y)
\end{array}\right)ds\ ,
\een where \ben\label{5.4a} F_1^{(1)}(z)=\mu_0 \lambda_0 e^{-\nu
_0z},\quad F_2^{(1)}(z)=\mu _0e^{-\nu _0z}\ . \een

On the other hand, using the representations (\ref{2.6}) and
(\ref{2.9}), the contribution of the right hand side of (\ref{5.2})
to the integral reads
 \beno &&\int\left (\psi_1-X_1^++\f ba\psi_2\right )\, e^{i\zeta y}\f
 {d\lambda}{2\pi \zeta }
=-\int_x^{+\infty}\Psi(x,s)\int X_1^+(s,\lambda)e^{i\zeta y}\f {d\lambda }{2\pi \zeta }\, ds\nonumber\\
&&\qquad+\int e^{i\zeta y} X_2^+(x,\lambda)\f b a(\lambda)\f
{d\lambda }{2\pi \zeta }-\int_x^{+\infty}\Psi(x,s)\int e^{i\zeta y}
X_2^+(s,\lambda)\f b a(\lambda)\f {d\lambda}{2\pi \zeta }\,  ds.
\eeno Note that, if $f$ is holomorphic and bounded on
$\overline{\Gamma }^+$,
 \ben\label{5.5} \int e^{i\zeta y}
f(\lambda)\f {d\lambda}{2\pi \zeta }= \int_{-\infty}^{+\infty}
e^{i\zeta y} \f {f(\lambda(\zeta ))-f(-\lambda(\zeta ))}{\lambda
(\zeta )}\f {d\zeta}{2\pi }\ . \een Then, we get by using
(\ref{5.5}) that \ben\label{5.6} &&\int e^{i\zeta y}\left
(\psi_1-X_1^++\f ba\psi_2\right )\f {d\lambda}{2\pi \zeta }
=-\Psi(x,y)\left(\begin{array}{c}
0\\
2\sqrt{2}
\end{array}\right)\nonumber\\
&&\qquad\qquad+\left(\begin{array}{c}
\sqrt{2}F_1^{(2)}(x+y)+i\sqrt{2}F_2^{(2)'}(x+y)\\
F_2^{(2)}(x+y)
\end{array}\right)\nonumber\\
&&\qquad\qquad-\int_x^\infty\Psi(x,s)\left(\begin{array}{c}
\sqrt{2}F_1^{(2)}(s+y)+i\sqrt{2}F_2^{(2)'}(s+y)\\
F_2^{(2)}(s+y)
\end{array}\right)ds\ ,
\een where \ben\label{5.6a}
F_1^{(2)}(z)=\f
1{2\pi}\int_{-\infty}^{+\infty}c_1(\zeta)e^{i\zeta z}d\zeta \ ,
F_2^{(2)}(z)=\f 1{2\pi}\int_{-\infty}^{+\infty}c_2(\zeta)e^{i\zeta
z}d\zeta \ . \een
 Combining (\ref{5.4}) and (\ref{5.6}), and
recalling the symmetry properties (\ref{symmetry}) of $\Psi $, this
completes the proof of Proposition \ref{marchenko}.
\endproof
\noindent Marchenko equations stated in this proposition are
related to the asymptotic analysis of Jost solutions for
$x\rightarrow +\infty $ and for this reason are called Marchenko
equations from the right. Marchenko equations from the left can be
derived similarly. We state the result without proof. Set
$$\widetilde c_1(\zeta)= \widetilde c(\lambda)+\widetilde c(-\lambda),\qquad \widetilde c_2(\zeta)=\f {\widetilde c(\lambda)-
\widetilde c(-\lambda)}{\lambda},\quad \widetilde c=-\f {b^*} a ,$$
and define the following real-valued functions of the real variable
$z$,
 \ben \label{5.11a}
  &&\widetilde F_1^{(1)}(z)= \mu _0\lambda _0e^{\nu_0z}\ ,\  \widetilde F_2^{(1)}(z)=\mu_0
 e^{\nu_0z}\ ,\\ \label{5.11b}
 &&\widetilde F_1^{2}(z)=\f
1{2\pi}\int_{-\infty}^{+\infty}\widetilde c_1(\zeta)e^{i\zeta
z}d\zeta \ ,\ \widetilde F_2^{(2)}(z)=\f
1{2\pi}\int_{-\infty}^{+\infty}\widetilde c_2(\zeta)e^{i\zeta
z}d\zeta \ ,\\ \label{5.11c} && \widetilde{F}_1=\widetilde{F}_1^{(2)}-\widetilde{F}_1^{(1)}\ ,\
\widetilde{F}_2=\widetilde{F}_2^{(2)}-\widetilde{F}_2^{(1)}\ .\een
\begin{Proposition}\label{marchenkoleft}
The functions $\Phi _{11}$, $\Phi _{12}$ defined by (\ref{2.11}),
(\ref{2.12}) satisfy the following system of equations,
\ben\label{5.10}
&&2\sqrt{2}{\Phi}_{11}(x,y)=-\widetilde{F}_2(x+y)\nonumber\\
&&\quad-\int^x_{-\infty}{\Phi}_{12}(x,s)\sqrt{2}(\widetilde{F}_1(s+y)-i\widetilde{F}_2'(s+y))+{\Phi}_{11}(x,s)\widetilde{F}_2(s+y)ds\
, \een and \ben\label{5.11}
&&2\sqrt{2}{\Phi}_{12}(x,y)=\sqrt{2}(\widetilde{F}_1(x+y)+i\widetilde{F}_2'(x+y))\nonumber\\
&&\quad-\int^x_{-\infty}{\Phi}_{11}(x,s)\sqrt{2}(\widetilde{F}_1(s+y)+i\widetilde{F}_2'(s+y))+{\Phi}_{12}(x,s)\widetilde{F}_2(s+y)ds\
, \een
for $y\le x$.
\end{Proposition}

\vskip 0.25cm \no In particular, in the unperturbed case case when
$q(x)=\f {\sqrt{2}} 2\tanh(\f x {\sqrt{2}})$, we have
\ben\label{5.12} a(\lambda)=\f {\lambda+\zeta-\f {\sqrt{2}} 2i}
{\lambda+\zeta+\f {\sqrt{2}} 2i}, \quad b(\lambda)=0,\quad
\lambda_0=0, \quad \mu_0=-2. \een Thus, $F_1(z)=0,\,\,
F_2(z)=2e^{-\f{\sqrt{2}}2z}$. The Marchenko equations become
\ben\label{5.13}
&&2\sqrt{2}\Psi_{11}(x,y)=2e^{-\f{\sqrt{2}}2(x+y)}\nonumber\\
&&\quad-2i\int_x^{+\infty}\Psi_{12}(x,s)e^{-\f{\sqrt{2}}2(y+s)}ds-
2\int_x^{+\infty}\Psi_{11}(x,s)e^{-\f{\sqrt{2}}2(y+s)}ds,
\een
and
\ben\label{5.14}
&&2\sqrt{2}\Psi_{12}(x,y)=-2ie^{-\f{\sqrt{2}}2(x+y)}\nonumber\\
&&\quad+2i\int_x^{+\infty}\Psi_{11}(x,s)e^{-\f{\sqrt{2}}2(y+s)}ds-
2\int_x^{+\infty}\Psi_{12}(x,s)e^{-\f{\sqrt{2}}2(y+s)}ds\ ,\een
which allows to recover the values of $\Psi ^0(x,y)$ already found
in subsection 2.3. (see (\ref{Psi0})).
\section{Further properties of the transition coefficients}

In this section, we will use the Marchenko equations to obtain the further information about the transition coefficients
which will play an important role in the proof of Theorem 1.1. In section 3, we proved that $a(\lambda)$ has at most one zero.
In the following, we will prove the existence of zero.

\begin{Proposition}\label{prop6.1}
The function $a(\lambda)$ has exactly one zero.
\end{Proposition}

\no Proof.\, By Lemma \ref{Lem3.2}, we know that $a(\lambda)$ has at
most one zero $\lambda_0$, and  that $|\lambda_0|\le C\vep$. It
remains to prove the existence of this zero. Assume that
$a(\lambda)$ does not vanish. In this case, the Marchenko equations
become \beno
&&2\sqrt{2}\Psi_{11}(x,y)=F_2^{(2)}(x+y)\nonumber\\
&&\quad-\int_x^{+\infty}\Psi_{12}(x,s)\sqrt{2}(F_1^{(2)}(s+y)-iF_2^{(2)'}(s+y))+\Psi_{11}(x,s)F_2^{(2)}(s+y)ds,
\eeno
and
\beno
&&2\sqrt{2}\Psi_{12}(x,y)=\sqrt{2}(F_1^{(2)}(x+y)+iF_2^{(2)'}(x+y))\nonumber\\
&&\quad-\int_x^{+\infty}\Psi_{11}(x,s)\sqrt{2}(F_1^{(2)}(s+y)+iF_2^{(2)'}(s+y))+\Psi_{12}(x,s)F_2^{(2)}(s+y)ds.
\eeno Thus, we get by the Young inequality that \ben\label{6.1}
&&\|\Psi(0,\cdot)\|_{L^2(y\ge 0)}\le C(1+\|\Psi (0,.)\|_{L^1(y\ge
0)}) (\|F_1^{(2)}\|_{L^2}+\|F_2^{(2)}\|_{L^2}+\|F_2^{(2)'}\|_{L^2}).
\een We get by (\ref{Phi}) and  Lemma \ref{PsiPhiperturb} that
\ben\label{6.2} \|\Psi(0,.)\|_{L^1(y\ge 0)}\le C,\quad
\|\Psi(0,\cdot)-\Psi^0(0,\cdot)\|_{L^2(y\ge 0)}\le C\vep\ , \een
where $\Psi ^0$ denotes the unperturbed kernel, given by
(\ref{Psi0}).  By Lemma \ref{Lem3.1}, we have \ben\label{6.3}
\|F_1^{(2)}\|_{L^2}+\|F_2^{(2)}\|_{L^2}+\|F_2^{(2)'}\|_{L^2}\le
C\vep^\f12. \een Indeed, by the Plancherel formula, \beno
\|F_1^{(2)}\|_{L^2}^2\le C\int_{\R}\left |\f {b(\lambda)}
{a(\lambda)}\right |^2d\zeta =C\int_{|\zeta|\le \varepsilon }d\zeta
+C\vep^2\int_{|\zeta|\ge \varepsilon }\f 1 {\zeta^2}d\zeta \le
C\vep. \eeno By summing up (\ref{6.1})-(\ref{6.3}), we obtain
$$
\f {\sqrt{2}} 4=\|\Psi^0(0,\cdot)\|_{L^2(y\ge 0)}^2\le C\vep,
$$
which is impossible. So, $a(\lambda)$ must have a zero.\endproof

\begin{Proposition}\label{prop6.2}
If $\vep$ is small enough, then the transition coefficients
$a(\lambda),b(\lambda)$ satisfy, on $H$,
 \ben\label{6.4}
\lim_{\zeta\rightarrow 0}\zeta a(\lambda)=0,\quad
\lim_{\zeta\rightarrow 0}\zeta b(\lambda)=0. \een
\end{Proposition}

In order to prove Proposition \ref{prop6.2}, we need the following two lemmas.
\begin{Lemma}\label{Lem6.1}
Let $\{F_1^{(1)}(x), F_2^{(1)}(x)\}$ be given by (\ref{5.4a}),
and $\{\widetilde{F}_1^{(1)}(x), \widetilde{F}_2^{(1)}(x)\}$ be given by (\ref{5.11a}). Then there hold
\beno
&&|F_1^{(1)}(x)|+|F_2^{(1)}(x)+2e^{-\f{\sqrt{2}}2x}|+|F_2^{(1)'}(x)-\sqrt{2}e^{-\f{\sqrt{2}}2x}|
\le C\vep e^{-\f 12x}\quad \textrm{for}\,\, x\ge 0,\\
&&|\widetilde{F}_1^{(1)}(x)|+|\widetilde{F}_2^{(1)}(x)+2e^{\f{\sqrt{2}}2x}|+|\widetilde{F}_2^{(1)'}(x)+\sqrt{2}e^{\f{\sqrt{2}}2x}|
\le C\vep e^{\f 12x}\quad \textrm{for}\,\, x\le 0.
\eeno
\end{Lemma}

\no Proof.\,\,Lemma \ref{Lem6.1} is a simple consequence of Lemma
\ref{Lem3.2}.\endproof

\begin{Lemma} \label{Lem6.2}
Let $\{F_1^{(2)}(x), F_2^{(2)}(x)\}$ be given by (\ref{5.6a}), and
$\{\widetilde{F}_1^{(2)}(x), \widetilde{F}_2^{(2)}(x)\}$ be given by
(\ref{5.11b}). Then there exists $M>0$ (independent of $\vep $ ) such that \ben
&&\int _M^{+\infty }\Bigl(|F_1^{(2)}(x)|+|F_2^{(2)}(x)|+|F_2^{(2)'}(x)|\Bigr)\, dx\le C\vep \quad ,\label{6.5}\\
&&\int _{-\infty }^{-M
}\Bigl(|\widetilde{F}_1^{(2)}(x)|+|\widetilde{F}_2^{(2)}(x)|+|\widetilde{F}_2^{(2)'}(x)|\Bigr)\,
dx\le C\vep \ .\label{6.6} \een
\end{Lemma}

\no Proof.\, We just prove (\ref{6.5}), since (\ref{6.6}) can be
similarly proved. We use the notation of Lemma \ref{PsiPhiperturb}.
Substracting the Marchenko equations for $\Psi ^0$ from the
Marchenko equations for $\Psi $, we obtain \beno
&&2\vep \sqrt{2}\Psi_{11}^1(x,y)=F_2^{(2)}(x+y)-(F^{(1)}_2(x+y)+2e^{-\f{\sqrt{2}}2(x+y)})\nonumber\\
&&\quad-\int_x^{+\infty}\Psi_{12}(x,s)\sqrt{2}(F_1^{(2)}(s+y)-iF_2^{(2)'}(s+y))+\Psi_{11}(x,s)F_2^{(2)}(s+y)ds\nonumber\\
&&\quad+\int_x^{+\infty}\Psi_{12}(x,s)[F_1^{(1)}(s+y)+i(2e^{-\f{\sqrt{2}}2(y+s)}-\sqrt{2}F_2^{(1)'}(s+y))]ds\nonumber\\
&&\quad+\int_x^{+\infty}\Psi_{11}(x,s)(2e^{-\f{\sqrt{2}}2(y+s)}+F_2^{(1)}(s+y))ds\nonumber\\
&&\quad-\int_x^{+\infty}2i\vep
\Psi_{12}^1(x,s)e^{-\f{\sqrt{2}}2(y+s)}+2\vep
\Psi_{11}^1(x,s)e^{-\f{\sqrt{2}}2(y+s)}ds, \eeno and \beno &&2\vep
\sqrt{2}\Psi_{12}^1(x,y)=\sqrt{2}(F_1^{(2)}+iF_2^{(2)'})(x+y)
-\sqrt{2}[F_1^{(1)}(x+y)+i(F_2^{(1)'}(x+y)-\sqrt{2}e^{-\f{\sqrt{2}}2(x+y)})]\nonumber\\
&&\quad-\int_x^{+\infty}\Psi_{11}(x,s)\sqrt{2}(F_1^{(2)}(s+y)+iF_2^{(2)'}(s+y))+\Psi_{12}(x,s)F_2^{(2)}(s+y)ds\nonumber\\
&&\quad+\int_x^{+\infty}\Psi_{11}(x,s)[F_1^{(1)}-i(2e^{-\f{\sqrt{2}}2(y+s)}-\sqrt{2}F_2^{(1)'}(s+y))]ds\nonumber\\
&&\quad+\int_x^{+\infty}\Psi_{12}(x,s)(2e^{-\f{\sqrt{2}}2(y+s)}+F_2^{(1)}(s+y))ds\nonumber\\
&&\quad+\int_x^{+\infty}2i\vep
\Psi_{11}^1(x,s)e^{-\f{\sqrt{2}}2(y+s)}-2\vep
\Psi_{12}^1(x,s)e^{-\f{\sqrt{2}}2(y+s)}ds. \eeno Integrating in
$y\in [x,+\infty [$, and using Lemma \ref{PsiPhiperturb} and Lemma \ref{Lem6.1}, we infer,
for every $x\ge 0$, \beno
&&\int _{2x}^{+\infty }\Bigl(|F_1^{(2)}(z)|+|F_2^{(2)}(z)|+|F_2^{(2)'}(z)|\Bigr)\, dz\le C\vep\\
&&\quad+C\int_0^{+\infty}|\underline \Psi (x,p)|\int _{2x+2p}^{+\infty
}\Bigl(|F_1^{(2)}(z)|+|F_2^{(2)}(z)|+|F_2^{(2)'}(z)|\Bigr)\, dz\,
dp. \eeno In view of estimate (\ref{Phi}), this completes the proof
of Lemma \ref{Lem6.2}\endproof

\vskip 0.25cm \noindent
Now we are in position to complete the proof of Proposition \ref{prop6.2}.
We can expand $a(\lambda)$ and $b(\lambda)$ near $\zeta=0$ as
\ben\label{6.7}
&&a(\lambda)=a^0(\zeta)+\f {\sigma_+(\varepsilon)} {\zeta}+a_1(\zeta),
\quad b(\lambda)=\f {\sigma_+(\vep)} {\zeta}+b_1(\zeta)\quad \textrm{for}\, \lambda>0,\\
&&a(\lambda)=a^0(\zeta)+\f {\sigma_-(\varepsilon)}
{\zeta}+a_1(\zeta), \quad b(\lambda)=-\f {\sigma_-(\vep)}
{\zeta}+b_1(\zeta)\quad  \textrm{for} \,\lambda<0.\label{6.8} \een
Here $ a^0(\zeta)=\f {\lambda+\zeta-\f {\sqrt{2}} 2i}
{\lambda+\zeta+\f {\sqrt{2}} 2i}, $ and from Lemma
\ref{psiphiperturb} and (\ref{5.1a}), we know that
$\sigma_\pm(\varepsilon)$ is a real constant and analytic in $\vep$
(tending to zero as $\vep\rightarrow 0$), $a_1(\zeta), b_1(\zeta)\in
C^1(\R)$ and for $k\le 1$ \ben |(\f d {d\zeta})^ka_1(\zeta)|+|(\f d
{d\zeta})^k b_1(\zeta)|\le C\vep.\label{6.9} \een From (\ref{6.7})
and (\ref{6.8}), we find that if $\sigma_\pm(\vep)\neq 0$, we have
\ben\label{6.10} \lim_{\lambda\rightarrow \f {\sqrt{2}} 2}\f
{b(\lambda)} {a(\lambda)}=1,\qquad \lim_{\lambda\rightarrow -\f
{\sqrt{2}} 2}\f {b(\lambda)} {a(\lambda)}=-1. \een

\begin{Lemma} \label{lem6.3}Let $\lambda=(\f12+\zeta^2)^\f12$.  There exists a constant
$C$ independent of $\delta_0$ and $\varepsilon$ such that as $\varepsilon$ tends to zero,
\ben
|\lim_{\delta\rightarrow 0}\f1{\pi}\int_{\delta\le |\zeta|\le \delta_0}\f {c(\lambda)} \zeta d\zeta
-i\textrm{sgn}\,\sigma_+(\varepsilon)|\le C\delta_0,\label{6.11}\\
|\lim_{\delta\rightarrow 0}\f1{\pi}\int_{\delta\le |\zeta|\le \delta_0}\f {c(-\lambda)} \zeta d\zeta
+i\textrm{sgn}\,\sigma_-(\varepsilon)|\le C\delta_0,\label{6.12}
\een
for $\delta_0>0$ small enough.
\end{Lemma}

\no{ Proof.}\, We first prove (\ref{6.11}).
We get by making an expansion for $a^0(\zeta)$ near $\zeta=0$ that
\beno
a(\lambda)=-i+\f {\sigma_+(\varepsilon)} {\zeta}+a_1(\zeta)+a_2(\zeta),
\eeno
where $a_2(\zeta)$ is smooth and $a_2(0)=0$. Thus, we get by (\ref{6.7}) that
\ben\label{6.15}
\lim_{\delta\rightarrow 0}\f1{\pi}\int_{\delta\le |\zeta|\le \delta_0}\f {c(\lambda)} \zeta d\zeta
&=&\lim_{\delta\rightarrow 0}\f1{\pi}\int_{\delta\le |\zeta|\le \delta_0}
\f {\f {\f {\sigma_+} {\zeta}+b_1(\zeta)} {-i+\f {\sigma_+(\varepsilon)} {\zeta}+a_1(\zeta)+a_2(\zeta)}} \zeta d\zeta\nonumber\\
&=&\lim_{\delta\rightarrow 0}\f1{\pi}\int_{\delta\le |\zeta|\le \delta_0}
\f {\f {\f {\sigma_+} {\zeta}+b_1(\zeta)} {-i+\f {\sigma_+(\varepsilon)} {\zeta}+a_1(\zeta)+a_2(\zeta)}-1} \zeta d\zeta\nonumber\\
&=& \lim_{\delta\rightarrow 0}\f1{\pi}\int_{\delta\le |\zeta|\le \delta_0}
\f {\f {i+b_1(\zeta)-a_1(\zeta)-a_2(\zeta)} {-i+\f {\sigma_+(\varepsilon)} {\zeta}+a_1(\zeta)+a_2(\zeta)}} \zeta d\zeta.
\een
Obviously, if $f(\zeta)=0(|\zeta|)$, then there holds for $\delta_0$ small enough
\beno
|\int_{\delta\le |\zeta|\le \delta_0}
\f {\f {f(\zeta)} {-i+\f {\sigma_+(\varepsilon)} {\zeta}+a_1(\zeta)+a_2(\zeta)}} \zeta d\zeta|\le C\delta_0,
\eeno
which together with (\ref{6.9})  implies that
\ben\label{6.16}
|\lim_{\delta\rightarrow 0}\f1{\pi}\int_{\delta\le |\zeta|\le \delta_0}\f {c(\lambda)} \zeta d\zeta
-\lim_{\delta\rightarrow 0}\f1{\pi}\int_{\delta\le |\zeta|\le \delta_0}
\f {\f {i+b_1(0)-a_1(0)} {-i+\f {\sigma_+(\varepsilon)} {\zeta}+a_1(\zeta)+a_2(\zeta)}} \zeta d\zeta|\le C\delta_0.
\een
Now, we write
\ben\label{6.17}
&&\lim_{\delta\rightarrow 0}\f1{\pi}\int_{\delta\le |\zeta|\le \delta_0}
\f {\f 1 {-i+\f {\sigma_+(\varepsilon)} {\zeta}+a_1(\zeta)+a_2(\zeta)}} \zeta d\zeta\nonumber\\
&&=\lim_{\delta\rightarrow 0}\f1{\pi}\int_{\delta\le |\zeta|\le \delta_0}
\f {\f {-a_1(\zeta)-a_2(\zeta)} {(-i+\f {\sigma_+(\varepsilon)} {\zeta}+a_1(\zeta)+a_2(\zeta))(-i+\f {\sigma_+(\varepsilon)} {\zeta})}} \zeta d\zeta
+\lim_{\delta\rightarrow 0}\f1{\pi}\int_{\delta\le |\zeta|\le \delta_0}
\f {\f 1 {-i+\f {\sigma_+(\varepsilon)} {\zeta}}} \zeta d\zeta.
\een
By a direct computation, we get for $\sigma_+(\varepsilon)\neq 0$
\ben
\lim_{\delta\rightarrow 0}\f1{\pi}\int_{\delta\le |\zeta|\le \delta_0}
\f {\f 1 {-i+\f {\sigma_+(\varepsilon)} {\zeta}}} \zeta d\zeta
&=&\lim_{\delta\rightarrow 0}\f1{\pi}\int_{\delta\le |\zeta|\le \delta_0}
\f {i\zeta+\sigma_+(\vep)} {\zeta^2+\sigma_+(\varepsilon)^2}  d\zeta\nonumber\\
&=&\lim_{\delta\rightarrow 0}\f2{\pi}\int_{\delta}^{\delta_0}
\f {\sigma_+(\vep)} {\zeta^2+\sigma_+(\varepsilon)^2} d\zeta\nonumber\\
&=& \f2{\pi}\int_{0}^{\f {\delta_0} {\sigma_+(\vep)}} \f 1
{1+\zeta^2} d\zeta=\f2{\pi}\arctan(\f {\delta_0} {\sigma_+(\vep)}),
\een and for $\sigma_+(\varepsilon)=0$ \ben \lim_{\delta\rightarrow
0}\f1{\pi}\int_{\delta\le |\zeta|\le \delta_0} \f {\f 1 {-i+\f
{\sigma_+(\varepsilon)} {\zeta}}} \zeta d\zeta=0. \een By the
properties of $a_1(\zeta)$ and $a_2(\zeta)$, it is easy to prove
that for $\sigma_+(\varepsilon)\neq 0$ \ben
&&|\lim_{\delta\rightarrow 0}\f1{\pi}\int_{\delta\le |\zeta|\le
\delta_0} \f {\f {-a_1(\zeta)-a_2(\zeta)} {(-i+\f
{\sigma_+(\varepsilon)} {\zeta}+a_1(\zeta)+a_2(\zeta))(-i+\f
{\sigma_+(\varepsilon)} {\zeta})}} \zeta d\zeta| \le
C\delta_0+C\varepsilon\ln \sigma_+(\vep), \een and for
$\sigma_+(\varepsilon)=0$ \ben\label{6.20}
&&|\lim_{\delta\rightarrow 0}\f1{\pi}\int_{\delta\le |\zeta|\le
\delta_0} \f {\f {-a_1(\zeta)-a_2(\zeta)} {(-i+\f
{\sigma_+(\varepsilon)} {\zeta}+a_1(\zeta)+a_2(\zeta))(-i+\f
{\sigma_+(\varepsilon)} {\zeta})}} \zeta d\zeta| \le C\delta_0. \een
Note that if $\sigma_+(\vep)\neq 0$, we have, by the analyticity of
$\sigma_+(\vep)$,
$$
\lim_{\vep\rightarrow 0}\varepsilon\ln \sigma_+(\vep)=0.
$$
By summing up (\ref{6.15})-(\ref{6.20}), we get that as $\vep$ tends to zero
\beno
|\lim_{\delta\rightarrow 0}\f1{\pi}\int_{\delta\le |\zeta|\le \delta_0}\f {c(\lambda)} \zeta d\zeta
-i\textrm{sgn}\,\sigma_+(\varepsilon)|\le C\delta_0,
\eeno
 for $\delta_0$ small enough. Since (\ref{6.12}) can be similarly proved, we omit its proof here.\endproof

To complete the proof of Proposition \ref{prop6.2}, it is sufficient to prove that $\sigma_\pm(\vep)$ must be zero.

\textbf{Case 1}. $\sigma_+(\vep)>0$.

By the definition of $F_2^{(2)}(x)$, we have
\ben\label{6.21}
\int_M^{+\infty}F_2^{(2)}(x)dx=\f12(\breve{c}_2(0)-i\textrm{H}\breve{c}_2(0)),
\een
where $\textrm{H}$ is the Hilbert transform and
\beno
\breve{c}_2(\zeta)=e^{iM\zeta}\f {c(\lambda)-c(-\lambda)}{\lambda},\quad
c=\f b a.
\eeno
We get by Lemma \ref{Lem6.2} that
\ben\label{6.22}
|\breve{c}_2(0)-i\textrm{H}\breve{c}_2(0)|\le C\vep.
\een
By the definition of the Hilbert transform, we have
\ben\label{6.23}
\textrm{H}\breve{c}_2(0)=\lim_{\delta\rightarrow 0}\f1{\pi}\int_{|\zeta|\ge \delta}\f {\breve{c}_2(\zeta)} \zeta d\zeta
=\lim_{\delta\rightarrow 0}\f1{\pi}\int_{\delta\le |\zeta|\le \delta_0}\f {\breve{c}_2(\zeta)} \zeta d\zeta
+\f1{\pi}\int_{|\zeta|\ge \delta_0}\f {\breve{c}_2(\zeta)} \zeta d\zeta,
\een
where $\delta_0$ is a small enough constant. We get by (\ref{3.8}) that
\ben\label{6.24}
|\f1{\pi}\int_{|\zeta|\ge \delta_0}\f {\breve{c}_2(\zeta)} \zeta d\zeta|\le C\vep/\delta_0.
\een
We write
\beno
\lim_{\delta\rightarrow 0}\f1{\pi}\int_{\delta\le |\zeta|\le \delta_0}\f {\breve{c}_2(\zeta)} \zeta d\zeta
&=&\lim_{\delta\rightarrow 0}\f1{\pi}\int_{\delta\le |\zeta|\le \delta_0}(\f {e^{iM\zeta}}{\lambda}-\sqrt{2})\f {c(\lambda)-c(-\lambda)} \zeta d\zeta\\
&&+\sqrt{2}\lim_{\delta\rightarrow 0}\f1{\pi}\int_{\delta\le |\zeta|\le \delta_0}\f {c(\lambda)-c(-\lambda)} \zeta d\zeta.
\eeno
Since $|c(\lambda)|\le 1$ by (\ref{3.8}), we get
\ben\
|\lim_{\delta\rightarrow 0}\f1{\pi}\int_{\delta\le |\zeta|\le \delta_0}(\f {e^{iM\zeta}}{\lambda}-\sqrt{2})\f {c(\lambda)-c(-\lambda)} \zeta d\zeta|\le C\delta_0.
\een
 From Lemma \ref{lem6.3}, we find that as $\varepsilon\rightarrow 0$
\ben\label{6.26}
|\lim_{\delta\rightarrow 0}\f1{\pi}\int_{\delta\le |\zeta|\le \delta_0}\f {c(\lambda)-c(-\lambda)} \zeta d\zeta
-i(1+\textrm{sgn}\,\sigma_-(\varepsilon))|\le C\delta_0.
\een

By summing up (\ref{6.23})-(\ref{6.26}), we get that as $\varepsilon\rightarrow 0$
\ben\label{6.27}
|\textrm{H}\breve{c}_2(0)-i\sqrt{2}(1+\textrm{sgn}\,\sigma_-(\varepsilon))|=o(1).
\een
On the other hand, we get by (\ref{6.8}) and (\ref{6.10}) that as $\varepsilon\rightarrow 0$
$$|\breve{c}_2(0)-\sqrt{2}(1+|\textrm{sgn}\,\sigma_-(\varepsilon)|)|=o(1),$$
from which and (\ref{6.27}), it follows that, as $\varepsilon\rightarrow 0$,
\beno
|\breve{c}_2(0)-i\textrm{H}\breve{c}_2(0)|=|2\sqrt{2}+\sqrt{2}(|\textrm{sgn}\,\sigma_-(\varepsilon)|+\textrm{sgn}\,\sigma_-(\varepsilon))|+o(1)\ge \sqrt{2},
\eeno
which contradicts (\ref{6.22}).

\textbf{Case 2}. $\sigma_-(\vep)>0$.

Exactly as in the proof of Case 1, we can get, as $\varepsilon\rightarrow 0$,
\beno
|\breve{c}_2(0)-i\textrm{H}\breve{c}_2(0)|
=|2\sqrt{2}+\sqrt{2}(|\textrm{sgn}\,\sigma_+(\varepsilon)|+\textrm{sgn}\,\sigma_+(\varepsilon))|+o(1)\ge \sqrt{2},
\eeno
which contradicts (\ref{6.22}).

\textbf{Case 3}. $\sigma_+(\vep)<0$.

By the definition of $\widetilde{F}_2^{(2)}(x)$, we have
\beno
\int^{-M}_{-\infty}\widetilde{F}_2^{(2)}(x)dx=\f12(\breve{\widetilde{c}}_2(0)+i\textrm{H}\breve{\widetilde{c}}_2(0)),
\eeno
where $\textrm{H}$ is the Hilbert transform and
\beno
\breve{\widetilde{c}}_2(\zeta)=e^{-iM\zeta}\f {\widetilde{c}(\lambda)-\widetilde{c}(-\lambda)}{\lambda},\quad
\widetilde{c}=-\f {b^*} a.
\eeno
We get by (\ref{6.6}) that
\ben\label{6.28}
|\breve{\widetilde{c}}_2(0)+i\textrm{H}\breve{\widetilde{c}}_2(0)|\le C\vep.
\een
Exactly as in the proof of Case 1, we can get, as $\varepsilon\rightarrow 0$,
\beno
&&|\breve{\widetilde{c}}_2(0)+\sqrt{2}(1+|\textrm{sgn}\,\sigma_-(\varepsilon)|)|=o(1),\\
&&|\textrm{H}\breve{\widetilde{c}}_2(0)-i\sqrt{2}(1-\textrm{sgn}\,\sigma_-(\varepsilon))|=o(1).
\eeno
So, we have
\beno
|\breve{\widetilde{c}}_2(0)+i\textrm{H}\breve{\widetilde{c}}_2(0)|=|2\sqrt{2}
+\sqrt{2}(|\textrm{sgn}\,\sigma_-(\varepsilon)|-\textrm{sgn}\,\sigma_-(\varepsilon))|+o(1)\ge \sqrt{2},
\eeno
which contradicts (\ref{6.28}).

\textbf{Case 4}. $\sigma_-(\vep)<0$.

As in case 3, we have, as $\varepsilon\rightarrow 0$,
\beno
|\breve{\widetilde{c}}_2(0)+i\textrm{H}\breve{\widetilde{c}}_2(0)|=|2\sqrt{2}
+\sqrt{2}(|\textrm{sgn}\,\sigma_+(\varepsilon)|-\textrm{sgn}\,\sigma_+(\varepsilon))|+o(1)\ge \sqrt{2},
\eeno
which contradicts (\ref{6.28}).

So, we conclude that $\sigma_\pm$ must be zero. This completes the proof of Proposition \ref{prop6.2}. \endproof
\vspace{0.2cm}

As a corollary of Proposition \ref{prop6.2}, we have
\begin{Proposition}\label{prop6.3} For $\vep$ small enough, we have
\ben\label{6.24} |\f {b(\lambda)} {a(\lambda)}|\le C\vep\langle
\zeta\rangle ^{-3}. \een
\end{Proposition}

\no Proof.\,By Lemma \ref{Lem3.1}, it suffices to prove that
\beno
|\f {b(\lambda)} {a(\lambda)}|\le C\vep,
\eeno
which can be deduced from (\ref{6.7}) and (\ref{6.8}), since $\sigma_\pm(\vep)$ is zero by Proposition \ref{prop6.2}.\endproof

\section{Proof of Theorem 1.1}
Let $\{a(\lambda), b(\lambda)\}$ be given by (\ref{3.3}). By Proposition \ref{prop6.1}, $a(\lambda)$ has only one zero
$\lambda_0$. Let $b_0$ be the constant such that
\beno
\varphi_1(x,\lambda_0)=b_0\psi_2(x,\lambda_0).\label{5.1}
\eeno
The evolution transition coefficients $\{a(t,\lambda), b(t,\lambda)\}$ are given by
\beno\label{7.1}
a(t,\lambda)=a(\lambda), \quad  b(t,\lambda)=b(\lambda)\exp(-i4\lambda\zeta t).
\eeno
and
\beno\label{7.2}
\mu_0(t)=\f {b_0(t)} {\nu_0a'(\lambda_0)},\quad b_0(t)=b_0\exp(4\lambda_0\nu_0 t),
\quad \zeta_0=(\lambda_0^2-\f12)^\f12=i\nu_0.
\eeno
We set
\begin{align}
&F_1^{(1)}(t,z)=\mu_0(t) \lambda_0 e^{-\nu_0z},\quad  F_2^{(1)}(t,z)=\mu_0(t) e^{-\nu_0z},\nonumber\\
&F_1^{(2)}(t,z)=\f 1{2\pi}\int_{-\infty}^{+\infty}c_1(t,\zeta)e^{i\zeta z}d\zeta,\quad
F_2^{(2)}(t,z)=\f 1{2\pi}\int_{-\infty}^{+\infty}c_2(t,\zeta)e^{i\zeta z}d\zeta,\nonumber\\
&c_1(t,\zeta)= c(t,\lambda)+c(t,-\lambda),\qquad c_2(t,\zeta)=\f {c(t,\lambda)-c(t,-\lambda)}{\lambda},\quad
c(t,\lambda)=\f {b(t,\lambda)} {a(t,\lambda)},\nonumber\\
&F_1(t,z)=F_1^{(2)}(t,z)-F_1^{(1)}(t,z),\qquad F_2(z)=F_2^{(2)}(t,z)-F_2^{(1)}(t,z).\nonumber
\end{align}
In the following, we will solve the evolution Marchenko equations for $y\ge x$
\ben\label{7.1}
&&2\sqrt{2}\Psi_{11}(t,x,y)=F_2(t,x+y)\nonumber\\
&&\quad-\int_x^{+\infty}\Psi_{12}(t,x,s)\sqrt{2}(F_1(t,s+y)-iF_2'(t,s+y))+\Psi_{11}(t,x,s)F_2(t,s+y)ds,
\een
and
\ben\label{7.2}
&&2\sqrt{2}\Psi_{12}(t,x,y)=\sqrt{2}(F_1(t,x+y)+iF_2'(t,x+y))\nonumber\\
&&\quad-\int_x^{+\infty}\Psi_{11}(t,x,s)\sqrt{2}(F_1(t,s+y)+iF_2'(t,s+y))+\Psi_{12}(t,x,s)F_2(t,s+y)ds.
\een Let $\{\Psi_{11}^{(1)}(t,x,y),\Psi_{12}^{(1)}(t,x,y)\}$ be the
solution of the following Marchenko equations involving only kernel
$F^{(1)}$, namely \beno
&&2\sqrt{2}\Psi_{11}^{(1)}(t,x,y)=-F_2^{(1)}(t,x+y)\\
&&\quad+\int_x^{+\infty}\Psi_{12}^{(1)}(t,x,s)\sqrt{2}(F_1^{(1)}(t,s+y)-iF_2^{(1)'}(t,s+y))+\Psi_{11}^{(1)}(t,x,s)F_2^{(1)}(t,s+y)ds,
\eeno and \beno
&&2\sqrt{2}\Psi_{12}^{(1)}(t,x,y)=-\sqrt{2}(F_1^{(1)}(t,x+y)+iF_2^{(1)'}(t,x+y))\\
&&\quad+\int_x^{+\infty}\Psi_{11}^{(1)}(t,x,s)\sqrt{2}(F_1^{(1)}(t,s+y)+iF_2^{(1)'}(t,s+y))+\Psi_{12}^{(1)}(t,x,s)F_2^{(1)}(t,s+y)ds.
\eeno Notice that, in view of the expression of $F^{(1)}_1,
F^{(1)}_2$ recalled above, these integral are finite rank equations.
Consequently, they can be exactly solved as \ben\label{7.3}
\Psi_{11}^{(1)}(t,x,y)=\f {\nu_0e^{\nu_0(x-y)}}
{1-\f{2\sqrt{2}\nu_0} {\mu_0(t)}e^{2\nu_0 x}}, \quad
\Psi_{12}^{(1)}(t,x,y)=\f
{\sqrt{2}\nu_0(\lambda_0-i\nu_0)e^{\nu_0(x-y)}}
{1-\f{2\sqrt{2}\nu_0} {\mu_0(t)}e^{2\nu_0 x}}. \een
\begin{Proposition}\label{prop7.1}
There exists a unique solution $\{\Psi_{11}(t,x,y),\Psi_{12}(t,x,y)\}$
to the Marchenko equations (\ref{7.1}) and (\ref{7.2}) such that for $|\alpha|\le 2$
\ben
&&\p^\alpha\Psi_{11}(t,x,y),\quad \p^\alpha\Psi_{12}(t,x,y)\in L^2_y(y\ge x),\label{7.4}\\
&&\|\p_y^\alpha\bigl(\Psi(t,x,\cdot)-\Psi^{(1)}(t,x,\cdot)\bigr)\|_{L^2(y\ge
x)}\le C\vep.\label{7.5} \een Furthermore, if we set \ben
u(t,x)=2\sqrt{2}i\Psi_{21}(t,x,x)+1,\label{7.6} \een then $u(t,x)$
is a solution of the Gross-Pitaevskii equation (\ref{1.1}).
\end{Proposition}

\no Proof.\quad We set \beno
\Psi_{11}^r(t,x,y)=\Psi_{11}(t,x,y)-\Psi_{11}^{(1)}(t,x,y),\quad
\Psi_{12}^r(t,x,y)=\Psi_{12}(t,x,y)-\Psi_{12}^{(1)}(t,x,y). \eeno
Then the Marchenko equations (\ref{7.1}) and (\ref{7.2}) are
transformed into the equations in terms of $\{ \Psi_{11}^r(t,x,y),
\Psi_{12}^r(t,x,y)\}$ \ben\label{7.7}
&&2\sqrt{2}\Psi_{11}^r(t,x,y)=\int_x^{+\infty}\Psi_{12}^r(t,x,s)\sqrt{2}(F_1^{(1)}(t,s+y)-iF_2^{(1)'}(t,s+y))\nonumber\\
&&\qquad+\Psi_{11}^r(t,x,s)F_2^{(1)}(t,s+y)ds+{\cal F}_1(t,x,y),
\een
and
\ben\label{7.8}
&&2\sqrt{2}\Psi_{12}^r(t,x,y)=\int_x^{+\infty}\Psi_{11}^r(t,x,s)\sqrt{2}(F_1^{(1)}(t,s+y)+iF_2^{(1)'}(t,s+y))\nonumber\\
&&\qquad+\Psi_{12}^r(t,x,s)F_2^{(1)}(t,s+y)ds+{\cal F}_2(t,x,y).
\een
Here
\ben
&&{\cal F}_1(t,x,y)=-\int_x^{+\infty}(\Psi_{12}^{(1)}+\Psi_{12}^r)(t,x,s)\sqrt{2}(F_1^{(2)}(t,s+y)-iF_2^{(2)'}(t,s+y))\nonumber\\
&&\quad+(\Psi_{11}^{(1)}+\Psi_{11}^r)(t,x,s)F_2^{(2)}(t,s+y)ds+F_2^{(2)}(t,x+y),\label{7.9}\\
&&{\cal F}_2(t,x,y)=-\int_x^{+\infty}(\Psi_{11}^{(1)}+\Psi_{11}^r)(t,x,s)\sqrt{2}(F_1^{(2)}(t,s+y)+iF_2^{(2)'}(t,s+y))\nonumber\\
&&\quad+(\Psi_{12}^{(1)}+\Psi_{12}^r)(t,x,s)F_2^{(2)}(t,s+y)ds+\sqrt{2}(F_1^{(2)}(t,x+y)+iF_2^{(2)'}(t,x+y)).\label{7.10}
\een
For fixed $t,x\in \R$, if the source terms $\{{\cal F}_1(t,x,y),{\cal F}_1(t,x,y)\}\in L^2_y(y\ge x)$ are given,
we firstly show that the integral equations (\ref{7.7}) and (\ref{7.8}) have a unique solution
$\{\Psi_{11}^r(t,x,y),\Psi_{12}^r(t,x,y)\}\in L^2_y(y\ge x)$.
We reformulate the integral equations (\ref{7.7}) and (\ref{7.8}) as
\ben
(I+\Omega_x)\Psi^r(t,x,y)={\cal F}(t,x,y),\label{7.11}
\een
where
$$
\Omega_x\Psi^r(t,x,y)=\int_x^{+\infty}\Omega(t,s+y)\Psi^r(t,x,s)ds.
$$
As we already noticed, the kernel $\Omega(t,s)$ is finite rank.
Therefore it suffices to show that the homogeneous equation
$$
(I+\Omega_x)\Psi^r(t,x,y)=0
$$
has only a trivial solution in $L^2_y(y\ge x)$.
Indeed, multiplying by $\overline{\Psi}_{11}^r(t,x,y)$ on both sides of (\ref{7.7}),
then integrating the resulting equations with respect to $y$ on $(x,+\infty)$, we obtain
\ben\label{7.12}
&&\int_x^{+\infty}|\Psi_{11}^r(t,x,y)|^2dy=\mu_0(t)|\int_x^{+\infty}e^{-\nu_0y}\Psi_{11}^r(t,x,y)dy|^2
+\int_x^{+\infty}{\cal F}_1(t,x,y)\overline{\Psi}_{11}^r(t,x,y)dy\nonumber\\
&&\qquad+\sqrt{2}\mu_0(t)(\lambda_0+i\nu_0)\int_x^{+\infty}e^{-\nu_0s}\Psi_{12}^r(t,x,s)ds\int_x^{+\infty}e^{-\nu_0y}\overline{\Psi}_{11}^r(t,x,y)dy.
\een
Similarly, we have
\ben\label{7.13}
&&\int_x^{+\infty}|\Psi_{12}^r(t,x,y)|^2dy=\mu_0(t)|\int_x^{+\infty}e^{-\nu_0y}\Psi_{12}^r(t,x,y)dy|^2
+\int_x^{+\infty}{\cal F}_2(t,x,y)\overline{\Psi}_{12}^r(t,x,y)dy\nonumber\\
&&\qquad+\sqrt{2}\mu_0(t)(\lambda_0-i\nu_0)\int_x^{+\infty}e^{-\nu_0s}\Psi_{11}^r(t,x,s)ds\int_x^{+\infty}e^{-\nu_0y}\overline{\Psi}_{12}^r(t,x,y)dy.
\een
We get by summing up (\ref{7.12}) and (\ref{7.13}) that
\ben\label{7.14}
&&2\sqrt{2}\int_x^{+\infty}(|\Psi_{11}^r(t,x,y)|^2+|\Psi_{12}^r(t,x,y)|^2)dy\nonumber\\
&&=\int_x^{+\infty}{\cal F}_1(t,x,y)\overline{\Psi}_{11}^r(t,x,y)dy+\int_x^{+\infty}{\cal F}_2(t,x,y)\overline{\Psi}_{12}^r(t,x,y)dy\nonumber\\
&&\quad+\mu_0(t)|\int_x^{+\infty}e^{-\nu_0s}\Psi_{11}^r(t,x,y)dy|^2+\mu_0(t)|\int_x^{+\infty}e^{-\nu_0y}\Psi_{12}^r(t,x,y)dy|^2\nonumber\\
&&\quad+2\sqrt{2}\mu_0(t)\lambda_0{\cal R}\Bigl(\int_x^{+\infty}e^{-\nu_0s}\Psi_{12}^r(t,x,s)ds\int_x^{+\infty}e^{-\nu_0y}\overline{\Psi}_{11}^r(t,x,y)dy\Bigr)\nonumber\\
&&\quad-2\sqrt{2}\mu_0(t)\nu_0{\cal I}\Bigl(\int_x^{+\infty}e^{-\nu_0s}\Psi_{12}^r(t,x,s)ds\int_x^{+\infty}e^{-\nu_0y}\overline{\Psi}_{11}^r(t,x,y)dy\Bigr).
\een
Note that the last two terms on the right hand side of (\ref{7.14}) are less than
\beno
&&2\sqrt{2}\mu_0(t)(\lambda_0^2+\nu_0^2)^\f12\Bigl|\int_x^{+\infty}e^{-\nu_0s}\Psi_{12}^r(t,x,s)ds\int_x^{+\infty}e^{-\nu_0y}\overline{\Psi}_{11}^r(t,x,y)dy\Bigr|\\
&&\quad\le \mu_0(t)\Bigl(|\int_x^{+\infty}e^{-\nu_0y}\Psi_{11}^r(t,x,y)dy|^2+|\int_x^{+\infty}e^{-\nu_0y}\Psi_{12}^r(t,x,y)dy|^2\Bigr).
\eeno
Thus, we obtain
\ben\label{7.15}
\int_x^{+\infty}(|\Psi_{11}^r(t,x,y)|^2+|\Psi_{12}^r(t,x,y)|^2)dy
\le C\int_x^{+\infty}(|{\cal F}_1(t,x,y)|^2+|{\cal F}_2(t,x,y)|^2)dy.
\een
In particular, if ${\cal F}_1(t,x,y)={\cal F}_2(t,x,y)=0$, then
$$
\Psi_{11}^r(t,x,y)=\Psi_{12}^r(t,x,y)=0.
$$
That is, the homogenous equation has a only trivial solution. Notice
that we proved in fact that the above integral equation is coercive,
which gives another argument for existence and uniqueness of a
solution in $L^2(y\geq x)$. Furthermore, it can be proved that if
$\p^\alpha\{{\cal F}_1(t,x,y),{\cal F}_1(t,x,y)\}\in L^2_y(y\ge x)$,
then there holds \beno \p^\alpha\Psi_{11}^r(t,x,y),\quad
\p^\alpha\Psi_{12}^r(t,x,y)\in L^2_y(y\ge x). \eeno Now, if  ${\cal
F}_1(t,x,y)$ and ${\cal F}_2(t,x,y)$ are given by (\ref{7.9}) and
(\ref{7.10}), it is easily verified by using Proposition
\ref{prop6.3} that for $|\alpha|\le 2$ \beno \|\p^\alpha{\cal
F}(t,x,\cdot)\|_{L^2(y\ge x)}\le C\vep +C\vep\sum_{|\beta|\le
|\alpha|}\|\p^\beta\Psi(t,x,\cdot)\|_{L^2(y\ge x)}. \eeno With this
and (\ref{7.15}), it can be proved by using the fixed point theorem
that the integral equations (\ref{7.7}) and (\ref{7.8}) has a unique
solution $\{\Psi_{11}^r(t,x,y),\Psi_{12}^r(t,x,y)\in L^2_y(y\ge
x)\}$ and there hold for $|\alpha|\le 2$ \beno
&&\p^\alpha\Psi_{11}^r(t,x,y),\quad \p^\alpha\Psi_{12}^r(t,x,y)\in L^2_y(y\ge x),\\
&&\|\p_y^\alpha\Psi^r_{11}(t,x,\cdot)\|_{L^2(y\ge x)}+\|\p_y^\alpha\Psi^r_{12}(t,x,\cdot)\|_{L^2(y\ge x)}\le C\vep.
\eeno
This completes the proof of the first part of Proposition \ref{prop7.1}.

Next, we show that $u(t,x)$ given by (\ref{7.6}) is a solution of (\ref{1.1}).
To prove it, we need the following two Lemmas whose proof will be given in the appendix.
\begin{Lemma}\label{Lem7.1} Let
\beno
\psi_1(t,x,\lambda)=X_1^+(x,\lambda)-\int_x^{+\infty} \Psi(t,x,s)X_1^+(s,\lambda)ds.
\eeno
Then $\psi_1(t,x,\lambda)$ is the solutions of (\ref{2.3}) with $q=\f {\sqrt{2}} 2u(t,x)$.
\end{Lemma}

\begin{Lemma}\label{Lem7.2} Let $\psi_1(t,x,\lambda)$ be as in Lemma \ref{Lem7.1}.
Then there holds
\ben
\f{\p \chi} {\p t}+iB_u\chi=i\sqrt{3}(\f E 2-\zeta)^2\chi,\label{7.16}
\een
with $\chi_1=(\sqrt{3}-1)^\f12e^{i\f {Ex} 2}\psi_{11},
\chi_2=(\sqrt{3}+1)^\f12e^{i\f {Ex} 2}\psi_{12}, \lambda=\f {\sqrt{3}} 2E$.
\end{Lemma}

With Lemma \ref{Lem7.1} and  Lemma \ref{Lem7.2}, we get by repeating the argument of section 4 that
$$
\f d {dt}L_u=i[L_u, B_u],
$$
which implies that $u(t,x)$ is a solution of  (\ref{1.1}).\endproof

\vspace{0.2cm}
With Proposition \ref{prop7.1}, we are in a position to complete the proof of Theorem 1.1.\vspace{0.2cm}

\no {\bf Proof of Theorem 1.1.}\,\, We recall that we denoted by
$\Psi ^0$ the unperturbed kernel given by (\ref{Psi0}). In view of
(\ref{Psi0}) and (\ref{7.3}), it is easy to verify, by using Lemma
\ref{Lem3.2}, that \ben\label{7.17}
\|\Psi^{(1)}(t,x+2\lambda_0t,y+2\lambda_0t)-\Psi^{0}(x,y)\|_{L^\infty(y\geq
x)}\le C\vep. \een We get by Proposition \ref{prop7.1} and the
Sobolev imbedding that \beno
\|\Psi(t,x,y)-\Psi^{(1)}(t,x,y)\|_{L^\infty(y\ge x)}\le C\vep, \eeno
from which and (\ref{7.17}), it follows that
$$
\|\Psi(t,x+2\lambda_0t,x+2\lambda_0t)-\Psi^{0}(x,x)\|_{L^\infty}\le
C\vep.
$$
Note that
$$
U_0(x)=\tanh(\f x{\sqrt{2}})=2\sqrt{2}i\Psi_{21}^0(x,x)+1.
$$
Thus, we have
\beno
\|u(t,x+2\lambda_0t)-U_0(x)\|_{L^\infty}\le 2\sqrt{2}\|\Psi(t,x+2\lambda_0t,x+2\lambda_0t)-\Psi^0(x,x)\|_{L^\infty}\le C\vep.
\eeno
This completes the proof of Theorem 1.1.\endproof
\section{Appendix A}

In this Appendix, we will prove Lemma \ref{Lem7.1} and Lemma \ref{Lem7.2}.\vspace{0.2cm}

\no {\bf Proof of Lemma \ref{Lem7.1}.}\,\,
To show that $\psi_1(t,x,\lambda)$ is the solution of (\ref{2.3}),
it suffices to prove that $\Psi(t,x,y)$ satisfies the linear system (\ref{2.7}). For this, we put
\beno
&&C_1(t,x,y)=(\p_x+\p_y)\Psi_{11}(t,x,y)-i\Bigl(-\f{\sqrt{2}}2\Psi_{12}(t,x,y)+q^*(t,x)\Psi_{12}^*(t,x,y)\Bigr),\\
&&C_2(t,x,y)=(\p_x-\p_y)\Psi_{12}(t,x,y)-i\Bigl(-\f{\sqrt{2}}2\Psi_{11}(t,x,y)+q^*(t,x)\Psi_{11}^*(t,x,y)\Bigr).
\eeno
Then by a direct calculation, we find that $C_1(t,x,y)$ and $C_2(t,x,y)$ satisfy the homogenous Marchenko equations
\beno
&&2\sqrt{2}C_1(t,x,y)=-\int_x^{+\infty}C_2(t,x,s)\sqrt{2}(F_1(t,s+y)-iF_2'(t,s+y))+C_1(t,x,s)F_2(t,s+y)ds,\\
&&2\sqrt{2}C_2(t,x,y)=-\int_x^{+\infty}C_1(t,x,s)\sqrt{2}(F_1(t,s+y)+iF_2'(t,s+y))+C_2(t,x,s)F_2(t,s+y)ds.
\eeno
Since the homogenous equations have a only trivial solution, $C_1(t,x,y)=C_2(t,x,y)=0$ follows.
This proves that $\Psi(t,x,y)$ satisfies the linear system (\ref{2.7}).\endproof

\vspace{0.2cm}
\no {\bf Proof of Lemma \ref{Lem7.2}.}\,\, Set
$h=e^{i\zeta x}\psi_1(t,x,\lambda)$. Then (\ref{7.16}) is equivalent to
\ben\label{8.1}
\f {\p h}{\p t}&=&-i\left(\begin{array}{cc}
\f{|u|^2-1}{\sqrt{3}+1} & 0\\
0 & \f{|u|^2-1}{\sqrt{3}-1}
\end{array}\right)h+\left(\begin{array}{cc}
0 & \f {\sqrt{3}+1}{\sqrt{2}}u^*_x\\
-\f {\sqrt{3}-1}{\sqrt{2}} u_x & 0
\end{array}\right)h\nonumber\\
&&+(2\sqrt{3}\zeta-2\lambda)\p_xh+i\sqrt{3}\p_x^2h\triangleq g.
\een
Since $\psi_1(t,x,\lambda)$ is a solution of (\ref{2.3}), we have
\beno
\lambda\p_xh=\left(\begin{array}{cc}
0 & q^*_x\\
q_x & 0
\end{array}\right)h+\left(\begin{array}{cc}
0 &\, q^*\\
q &\, 0
\end{array}\right)\p_xh+\left(\begin{array}{cc}
1 & 0\\
0 & -1
\end{array}\right)(i\p_x^2h+\zeta\p_xh).
\eeno
Substituting this into $g$, we get
\ben\label{8.2}
g&=&-i\left(\begin{array}{cc}
\f{|u|^2-1}{\sqrt{3}+1} & 0\\
0 & \f{|u|^2-1}{\sqrt{3}-1}
\end{array}\right)h+\left(\begin{array}{cc}
0 & \f {\sqrt{3}+1}{\sqrt{2}}u^*_x\\
-\f {\sqrt{3}-1}{\sqrt{2}} u_x & 0
\end{array}\right)h-2\left(\begin{array}{cc}
0 &\, q^*\\
q &\, 0
\end{array}\right)\p_xh\nonumber\\
&&+2\left(\begin{array}{cc}
\sqrt{3}-1 &\, 0\\
0 &\, \sqrt{3}+1
\end{array}\right)\zeta\p_xh+i\left(\begin{array}{cc}
\sqrt{3}-2 &\, 0\\
0 &\, \sqrt{3}+2
\end{array}\right)\p_x^2h.
\een
Note that
$$
h=\left(\begin{array}{c}
1\\
\sqrt{2}(\lambda-\zeta)
\end{array}\right)-\int_0^{+\infty}\Psi(t,x,x+s)X_1^+(s,\lambda)ds.
$$
Substituting this into (\ref{8.2}), we get by integrating by parts and (\ref{2.7}) that
\beno
g=\int_0^{+\infty}G(t,x,y+s)X_1^+(s,\lambda),
\eeno
where
\beno
G(t,x,y)=\left(\begin{array}{cc}
i(\p_x+\p_y)^2\Psi_{11}+2q^*(\p_x+\p_y)\Psi_{21} &\, i(\p_x+\p_y)^2\Psi_{12}+2q^*(\p_x+\p_y)\Psi_{22} \\
-i(\p_x+\p_y)^2\Psi_{21}+2q(\p_x+\p_y)\Psi_{11} &\, -i(\p_x+\p_y)^2\Psi_{22}+2q(\p_x+\p_y)\Psi_{12}
\end{array}\right).
\eeno
Note that
\beno
\f {\p h}{\p t}=-\int_0^{+\infty}\p_t\Psi(t,x,x+s)X_1^+(s,\lambda)ds.
\eeno
Thus, to prove (\ref{8.1}), it suffices to prove that
\ben\label{8.3}
\p_t\Psi(t,x,y)=-G(t,x,y).
\een
By the definition of $F_1(t,z)$ and $F_1(t,z)$, we have
\beno
&&\f {\p F_1}{\p t}(t,z)=-2\f {\p F_2}{\p z}(t,z)+2\f {\p^3 F_2}{\p z^3}(t,z),\\
&&\f {\p F_2}{\p t}(t,z)=-4\f {\p F_1}{\p z}(t,z).
\eeno
Thus, we get by differentiating (\ref{7.1}) and (\ref{7.2}) with respect to $t$ that
\beno
&&2\sqrt{2}\p_t\Psi_{11}(t,x,y)=D_1(t,x,y)\\
&&-\int_x^{+\infty}\p_t\Psi_{12}(t,x,s)\sqrt{2}(F_1(t,s+y)-iF_2'(t,s+y))+\p_t\Psi_{11}(t,x,s)F_2(t,s+y)ds,\\
&&2\sqrt{2}\p_t\Psi_{12}(t,x,y)=D_2(t,x,y)\\
&&-\int_x^{+\infty}\p_t\Psi_{11}(t,x,s)\sqrt{2}(F_1(t,s+y)+iF_2'(t,s+y))+\p_t\Psi_{12}(t,x,s)F_2(t,s+y)ds,
\eeno
where
\beno
D_1(t,x,y)&=&-\int_x^{+\infty}\Psi_{12}(t,x,s)\sqrt{2}(4F_2^{'''}(t,s+y)-2F_2'(t,s+y)+4iF_1^{''}(t,s+y))ds\\
&&+4\int_x^{+\infty}\Psi_{11}(t,x,s)F_1'(t,s+y)ds-4F_1'(t,x+y),\\
D_2(t,x,y)&=& \sqrt{2}(4F_2^{'''}(t,x+y)-2F_2'(t,x+y)-4iF_1^{''}(t,x+y))\\
&&-\int_x^{+\infty}\Psi_{11}(t,x,s)\sqrt{2}(4F_2^{'''}(t,s+y)-2F_2'(t,s+y)-4iF_1^{''}(t,s+y))ds\\
&&+4\int_x^{+\infty}\Psi_{12}(t,x,s)F_1'(t,s+y)ds.
\eeno
On the other hand, differentiating (\ref{7.1}) and (\ref{7.2}) with respect to $x$ and $y$,
we get by integrating by parts and (\ref{2.7}) that
\beno
&&2\sqrt{2}G_{11}(t,x,y)=-D_1(t,x,y)\nonumber\\
&&-\int_x^{+\infty}G_{12}(t,x,s)\sqrt{2}(F_1(t,s+y)-iF_2'(t,s+y))+G_{11}(t,x,s)F_2(t,s+y)ds,\\
&&2\sqrt{2}G_{12}(t,x,y)=-D_2(t,x,y)\nonumber\\
&&-\int_x^{+\infty}G_{11}(t,x,s)\sqrt{2}(F_1(t,s+y)+iF_2'(t,s+y))+G_{12}(t,x,s)F_2(t,s+y)ds.
\eeno
Thus, $\p_t\Psi(t,x,y)$ and $-G(t,x,y)$ satisfy the same Marchenko equations. So, we conclude
(\ref{8.3}) by the uniqueness of the solution of the Marchenko equations.\endproof

\section{Appendix B}

In this appendix, we prove Corollary \ref{cor}. The proof is based
on Theorem \ref{main} and a compactness argument.
\s

\no{\bf Proof.} Let $u_0$ be a Cauchy datum as in Theorem
\ref{main}. Due to Theorem \ref{main},  \beno \sup_{t\in
\R}\|u(x+2\lambda _0t,t)-U_0(x)\|_{L^\infty}\le C\varepsilon\ ,
\eeno where $\lambda _0$ is the unique zero of the transition
coefficient $a$ associated to $u_0$. We shall show that, given,
$\delta >0$, for $\varepsilon >0$ small enough, \beno \sup_{t\in
\R}\|u'(x+2\lambda _0t,t)-U_0'(x)\|_{L^2}+\||u(x+2\lambda
_0t,t)|^2-|U_0(x)|^2\|_{L^2}\le \delta. \eeno By contradiction,
otherwise  there would exist $\delta>0$ and a sequence $\{u_0^n\}$
verifying \beno \|<x>^4\p^k(u_0^n(x)-U_0(x))\|_{L^\infty}\rightarrow
0,\quad \textrm{as} \quad n\rightarrow+\infty,\qquad \textrm{for}
\quad k\le 3 \eeno and a sequence $\{t_n\}$ such that \ben
\label{9.1}
&&\|u_n'(t_n,x+2\lambda_0^nt_n)-U_0'(x)\|_{L^2}+\||u_n(t_n,x+2\lambda
_0t_n)|^2-|U_0(x)|^2\|_{L^2}\ge
\delta\ ,\\
&&\sup_{t\in
\R}\|u_n(t,x+2\lambda_0t)-U_0(x)\|_{L^\infty}\rightarrow 0,\quad
\textrm{as} \quad n\rightarrow+\infty.\label{9.2} \een Here $\lambda
_0^n$ denotes the unique zero of the transition coefficient $a_n$
associated to $u_0^n$, and $u_n$ denotes the solution to (\ref{1.1})
with Cauchy datum $u_0^n$.  We will follow the argument in Page 20
of \cite{Bet4} to yield a contradiction. Firstly, we have the energy
conservation
$$
E(u_n(t))=\f12\int|u_n'(t,x)|^2dx+\f14\int||u_n(t,x)|^2-1|^2dx=E(u_0^n),
$$
from which and (\ref{9.2}), it follows that,for any $B>0$,
\ben\label{9.3} u_n(t_n,x+2 \lambda_0^nt_n)\rightharpoonup
U_0(x)\quad \textrm{in}\quad  H^1(|x|\le B). \een Next, we take $B$
such that \ben\label{9.4} \int_{|x|\ge
B}\f12|U_0'|^2+\f14||U_0|^2-1|^2dx\le \f\delta 4, \een and
$n_\delta$ such that \ben\label{9.5} E(u_0^n)\le E(U_0)+\f\delta 4,
\een for any $n\ge n_\delta$. Since
$$
\int_{|x|\le B}\f12|U_0'|^2+\f14||U_0|^2-1|^2dx\le
\liminf_{n\rightarrow +\infty}\int_{-B+2 \lambda_0^nt_n}^{B+2\lambda
_0t_n}\f12|u_n'(t_n,x)|^2+\f14\|u_n(t_n,x)|^2-1|^2dx
$$
which together with (\ref{9.4})-(\ref{9.5}) and the energy
conservation implies that for $n\ge n_\delta$ \beno \int_{|x-2
\lambda_0^nt_n|\ge
B}\f12|u_n'(t_n,x)|^2+\f14\|u_n(t_n,x)|^2-1|^2dx\le \f\delta 2,
\eeno In view of (\ref{9.3}), we infer, for $n$ big enough \beno
\|u_n'(t_n,x+2 \lambda_0^nt_n)-U_0'(x)\|_{L^2}+\||u_n(t_n,x+2\lambda
_0t_n)|^2-|U_0(x)|^2\|_{L^2}<\delta \eeno which contradicts
(\ref{9.1}). The corollary follows.\endproof

\end{document}